\documentclass[review, 12pt]{elsarticle}

\usepackage{amssymb}
\usepackage{latexsym}
\usepackage{amsfonts}
\usepackage{amsmath}
\usepackage[latin1]{inputenc}
\usepackage{subfigure}
\usepackage{setspace}
\usepackage{graphicx}
\usepackage{graphics}
\usepackage{hyperref}
\usepackage{multirow}
\usepackage{color}
\usepackage{epsfig}
\usepackage{appendix}

\usepackage{soul} %this is to add commands \hl{} in order to 
\newcommand{\norm}[1]{\| #1 \|}
\newcommand{\eaa}[1]{\color{black}#1\normalcolor}

\journal{}

\begin{document}

\begin{frontmatter}



\title{A fourth-order exponential time differencing scheme with dimensional splitting for non-linear reaction-diffusion systems}


\author[label1]{E.O. Asante-Asamani\corref{cor1}}
\ead{easantea@clarkson.edu}

\address[label1]{Department of Mathematics, Clarkson University\\ Potsdam NY 13676, USA }

\author[label2,label3]{A. Kleefeld}
\ead{a.kleefeld@fz-juelich.de}

\address[label2]{J\"ulich Supercomputing Centre, Forschungszentrum J\"ulich GmbH,\\ 52425 J\"ulich, Germany}

\address[label3]{Faculty of Medical Engineering and Technomathematics, University of Applied Sciences Aachen,\\ 52428 J\"ulich, Germany}

\author[label4]{B.A. Wade}
\ead{bruce.wade@louisiana.edu}

\address[label4]{Department of Mathematics, University of Louisiana at Lafayette,\\ Lafayette LA 70504, USA }
\cortext[cor1]{Corresponding author}

\begin{abstract}
A fourth-order exponential time differencing (ETD) Runge-Kutta scheme with dimensional splitting is developed to solve multidimensional non-linear systems of reaction-diffusion equations (RDE). By approximating the matrix exponential in the scheme with the \eaa{A-acceptable} Pad\'e (2,2) rational function, the resulting scheme (ETDRK4P22-IF) is verified empirically to be fourth-order accurate for several RDE. The scheme is shown to be more efficient than competing fourth-order ETD and IMEX schemes, achieving up to 20X speed-up in CPU time. \eaa{Inclusion of up to three pre-smoothing steps of a lower order L-stable scheme facilitates efficient damping of spurious oscillations arising from problems with non-smooth initial/boundary conditions}. 

% Novel dimensional splitting techniques are developed for ETD Schemes which are second-order convergent and highly efficient. By using the ETD-Crank-Nicolson scheme we show that the proposed techniques can  reduce the computational time for non-linear reaction-diffusion systems by up to 70\%. Numerical tests are performed to empirically validate the superior performance of the splitting methods.
\end{abstract}


\end{frontmatter}


\section{Introduction}
Time-dependent reaction-diffusion equations (RDE) are mathematical models that describe the spatio-temporal dynamics of many physical processes in diverse fields such as cell and developmental biology \cite{mittal2021efficient,Kondo2010,landge2020pattern}, physics  \cite{ghafouri2020numerical}, geological engineering \cite{sundararajan2020groundwater}, and finance \cite{rambeerich2009exponential}. \eaa{Mathematically, RDE are described by the system of partial differential equations (PDE) of the form
\begin{equation}
\label{ADReq}
\begin{split}
\frac{\partial u}{\partial t} &= D\Delta u+ f(u,t), \quad \text{ in } \Omega_T \\
u(x,t) &=g(x,t), \text{ on } \Gamma_T.
\end{split}
\end{equation}
Here, $\Omega$ is a bounded open subset of $\mathbb{R}^d, d\in \mathbb{N}$, $\Omega_T:= \Omega\times (0,T]$ with parabolic boundary ${\Gamma_T:= \overline{\Omega_T}\setminus \Omega_T}$. The solution $u:\overline{\Omega_T}\rightarrow \mathbb{R}^s$, describes the concentration of $s$ species with diffusion coefficients $D = \mathrm{diag}(d_1,d_2,\ldots,d_s),\ d_i>0\ \forall i$. The function $f:\Omega_T\rightarrow \mathbb{R}^s$ describes the interaction between all the species and is assumed to be sufficiently smooth.} 

\eaa{In practice, the spatial derivatives in \eqref{ADReq} are first discretized on a partition of $\Omega$ with $m^d$ points to obtain a system of ordinary differential equations (ODE) 
\begin{align}
\label{ADReqmatrix}
\frac{dU}{dt} + AU &= F(U(t),t), \quad U\in \mathbb{R}^{s\times m^d}\\
U(x,0) &=U_0(x)\nonumber
\end{align} where $U(t) = (u_1(t),u_2(t),\ldots,u_s(t),\ldots,u(t)_{s\times m^d})^T$, $A$ is an $(s\times m^d)\times (s\times m^d)$ matrix approximation of $-D\Delta$ and $F\in \mathbb{R}^{s\times m^d}$ is typically a non-linear map approximating $f$ on the numerical grid for all $s$ species}. The major computational challenges in solving \eqref{ADReqmatrix} include numerical stiffness due to the spatial diffusion, coupling of very fast and very slow reaction kinetics, simulating problems in high spatial dimensions ($d\geq 2$) with non-linear source terms, spurious oscillations arising from non-smooth initial conditions or mismatched initial and boundary conditions as well as maintaining positive solutions especially for biological applications \cite{Hunsdorfer2003}.  

Among the host of time stepping methods available to solve stiff ODE systems \cite{Hunsdorfer2003,iserles2009first,shampine2018numerical,martin2019eserk5,martin2020solving} are the class of exponential time differencing (ETD) Runge-Kutta schemes.  These schemes approximate the exact solution of the system \eqref{ADReqmatrix},
\begin{equation}
\label{etdintegral}
     U(t) = \mathrm{e}^{-tA}U(0) + \int_0^t \mathrm{e}^{-(t-s)A}F(U(s),s) \;\mathrm{d}s\,,\quad t\in [0,T]
\end{equation} by approximating the integral and matrix exponential. In 2002, Cox and Matthews \cite{Cox2002} introduced a class of ETD Runge-Kutta (ETDRK) schemes which approximate the integral in \eqref{etdintegral} by interpolating the non-linear function $F$ with polynomials. The major attraction of these schemes is in the accurate treatment of the linear (diffusion) term via matrix exponential leading to stable evolution of stiff problems. 
 
 The computational time for ETDRK schemes has been improved by employing rational approximations of the matrix exponential. Several second-order and fourth-order ETDRK schemes have been developed that approximate the matrix exponential using Pad\'e rational functions (ETDRK-Pad\'e schemes) \cite{Kleefeld2012,Yousuf2012,yousuf2009efficient,Wade2009} as well as a rational function with real and distinct poles (RDP) leading to the ETDRDP scheme in \cite{asante2016real}. To further improve the computational efficiency of second-order ETDRK schemes when solving multidimensional RDE, we introduced a class of dimensionally split second-order ETDRK schemes to take advantage of fast solvers for tridiagonal linear systems. Utilizing Pad\'e (1,1) and RDP rational functions we developed efficient second-order ETDRK schemes with dimensional splitting in \cite{asante2016dimensional} and \cite{asante2020second}, respectively. The dimensional splitting scheme developed with RDP rational functions is referred to as ETDRDP-IF. \eaa{In spite of this progress, the savings in CPU time obtained from dimensional splitting has, until now, been limited to second-order ETDRK schemes}. 

 \eaa{In this work, we introduce a new class of fourth-order ETDRK schemes with dimensional splitting. By approximating the matrix exponential with the A-acceptable Pad\'e (2,2) rational function, we create an efficient fourth-order scheme we call ETDRK4P22-IF}.  Empirical convergence analysis validates the fourth-order accuracy of the scheme across various prototypical scalar and systems of RDE with Dirichlet and Neumann boundary conditions posed in $\mathbb{R}^2$. The new scheme is shown to significantly improve the run time of the existing fourth order ETDRK4P22 scheme (without dimensional splitting) developed by Yousuf \cite{Wade2009,yousuf2009efficient}, in some cases achieving 20X speed-up in CPU time with serial implementation of the algorithm. \eaa{The one-step nature of the scheme is also shown to provide a substantial computational advantage over competing fourth-order multistep IMEX schemes. By presmoothing ETDRK4P22-IF with a few steps of a lower order L-stable ETDRK scheme, we show that ETDRK4P22-IF can accurately solve RDE with non-smooth initial or boundary conditions.}  

\eaa{The remainder of the manuscript is organized as follows. In the Section \ref{sec:matandmethods} we describe the finite difference scheme used to obtain a fourth-order approximation the Laplacian. This is followed by details of the dimensional splitting of the fourth-order ETDRK scheme along with a rational approximation of the matrix exponential using Pad\'e(2,2) functions to yield the ETDRK4P22-IF scheme. Serial and parallel algorithms for implementing the ETDRK4P22-IF are also presented in this section. The convergence and efficiency of ETDRK4P22-IF is investigated in Section \ref{sec:experiments} for problems with and without exact solution. Finally, we wrap up with concluding remarks in Section \ref{sec:conclusion}.}

\section{Materials and Methods}
\label{sec:matandmethods}
\subsection{Discretization in space}
\label{sec:space}
To discretize the Laplacian $\Delta = \partial_{xx} + \partial_{yy}$ on the domain $[a,b]^2$, we partition the spatial domain in each direction into $m+2$ points with a uniform mesh of size $h=\frac{b-a}{m+1}$ with $m\geq 3.$ If we set $x_j=a + jh$ with $j=0,1,2,\ldots,m+1$, then the second-order partial derivative of a function $w(x,t)$ with respect to $x$ at $x_j$ can be discretized using the fourth-order central difference scheme (\cite{mathews2004numerical}, p. 339)
\begin{equation}
\label{dirichmain}
    \partial_{xx}w|_{x_j} = \frac{1}{12h^2}(-W_{j-2}+16W_{j-1}-30W_j + 16W_{j+1}-W_{j+2}) + \mathcal{O}(h^4)\,
\end{equation}
$\quad j=2,3,\ldots,m-1$, with $W_j\approx w(x_j,t)$.
The approximation used for the nodes at $j=0,1,m,m+1$ depends on the boundary condition. Since the solution at $x_0$ and $x_{m+1}$ are known for homogeneous Dirichlet boundary conditions, the derivatives at $x_1$ and $x_m$ are estimated with the formulas
\begin{equation}
\label{dirichmain1}
\begin{split}
    \partial_{xx}w|_{x_1} & \approx \frac{1}{12h^2}(11W_{0}-20W_1 + 6W_{2}+4W_3-W_{4})\,, \\
    \partial_{xx}w|_{x_m} & \approx \frac{1}{12h^2}(-W_{m-3}+4W_{m-2} + 6W_{m-1}-20W_{m}+11W_{m+1})\,   
\end{split}
\end{equation}
obtained by extrapolating the values, $W_{-1},W_{m+2}$ using a fourth degree Lagrange interpolating polynomial centered at $x_0$ and $x_{m+1}$.

% a fourth degree Lagrange interpolating polynomial centered at the end points $x_0$ and $x_{m+1}$ is used to extrapolate the values 
% \[W_{-1}=5W_0-10W_1+10W_2-5W_3+W_4\] 
% and 
% \[W_{m+2}=W_{m-3}-5W_{m-2}+10W_{m-1}-10W_m+5W_{m+1}\] 
% from which the
For homogeneous Neumann boundary conditions, we use the approximation
\begin{equation}
\label{nuemann1}
    \begin{split}
    \partial_{xx}w|_{x_1} & \approx \frac{1}{12h^2}(16W_{0}-31W_1 + 16W_{2}-W_3)\ ,\\
    \partial_{xx}w|_{x_m} & \approx \frac{1}{12h^2}(-W_{m-2}+16W_{m-1} - 30W_{m}+16W_{m+1}-W_{m+2}) \, 
\end{split}
\end{equation}
which is derived by setting \eaa{$W_{-1}= W_{1}$ and $W_m = W_{m+2}$}. The derivatives at $x_0$ and $x_{m+1}$ are set to 
\begin{equation}
\label{nuemann2}
\begin{split}
    \partial_{xx}w|_{x_0} &\approx \frac{1}{12h^2}(-30W_{0}+32W_1 -2W_{2})\,,\\
    \partial_{xx}w|_{x_{m+1}} &\approx \frac{1}{12h^2}(-2W_{m-1}+32W_{m} - 30W_{m+1})\,. 
\end{split}
\end{equation} Here the Neumann boundary condition is discretized using a fourth-order central difference approximation to the first-order partial derivative (see \cite{mathews2004numerical} p. 339) \eaa{with the additional assumption that $W_{-1}=W_{1}$ and $W_{m+2}=W_{m}$ }.  This discretization is shown empirically to be fourth-order accurate in our numerical experiments. A similar scheme has been used in \cite{gibou2005fourth} to achieve fourth-order accuracy in irregular domains. 
% is used to obtain $W_{-2} = W_2 - 8W_1+8W_{-1}$ which simplifies to $W_{-2} = W_2$ after substituting $W_{-1}=W_1$. Similarly, we obtain $W_{m+3}=W_{m-1}$. With these approximations the final set of partial derivatives are obtained as

The matrix approximating the second-order partial derivatives in the PDE in the presence of homogeneous Dirichlet boundary conditions has a bandwidth of six whereas the matrix for the case of Neumann boundary conditions has a bandwidth of five. To formulate the two dimensional Laplacian, we take advantage of the Kronecker product. Let $A=A_1 + A_2$ be the matrix approximation of the two dimensional (2D) Laplacian with $A_1$ and $A_2$ as the matrix approximations of $\partial_{xx}$ and $\partial_{yy}$ respectively in 2D. Let $B_p$ be the matrix representations of the second-order partial derivative in 1D with entries defined by equations \eqref{dirichmain}--\eqref{nuemann2} and $I_p$ an $p$-dimensional identity matrix. Then $A_1 = B_p \otimes I_p$ and  $A_2 = I_p \otimes B_p$ (see \cite{Hunsdorfer2003} for details). For Dirichlet boundary conditions and $s$ species, $p=s\cdotp m$, whereas $p=s\cdotp (m+2)$ for Neumann boundary.  We  proved in \cite{asante2020second,asante2016dimensional} that $A_1$ and $A_2$ commute. This is an important property which facilitates the formulation of the dimensional splitting scheme presented in the next section.

\subsection{Discretization in time}
\label{sec:time}
\subsubsection{Dimensional splitting of fourth-order ETD Runge-Kutta scheme }
In this section, we introduce a dimensional splitting version of the class of fourth-order ETD Runge-Kutta schemes developed by Cox-Matthews \cite{Cox2002} and Kassam-Trefethen \cite{kassam2005fourth} for solving the semi-linear system of ordinary differential equations in \eqref{ADReqmatrix}. Our motivation is to speed up the simulation time when solving problems in high spatial dimensions by taking advantage of a nicer matrix structure resulting from splitting dimensions.
% \begin{align}
% \label{odesystem}
%     u^{'}(t) +Au(t) &= F(u(t),t)\,,\quad  u(t)\in \mathbb{R}^m\,,\quad A \in \mathbb{R}^{m\times m}\,,\quad t\in [0,T]\\
%     u(0) &= u_0\,.\nonumber
% \end{align}
% In our applications, $A$ is a matrix resulting from the finite difference discretization of a spatial differential operator such as the Laplacian and $F: \R^m \rightarrow \R^m$ is a non-linear function.  The equation \eqref{odesystem} thus results from the application of the method of lines to some parabolic partial differential equation. 

The scheme is
\begin{equation}
\label{semidiscrete}
\begin{split}
    U_{n+1} &= \mathrm{e}^{-kA}U_n + \frac{1}{k^2}(-A)^{-3}[-4I+kA+\mathrm{e}^{-kA}(4I+3kA+k^2A^2)] F(U_n,t_n) \\
    &+   \frac{2}{k^2}(-A)^{-3}[2I-kA-\mathrm{e}^{-kA}(2+kA)]\left(F(a_n,t_n+\frac{k}{2})+F(b_n,t_n+\frac{k}{2})\right) \\
    &+ \frac{1}{k^2}(-A)^{-3}[-4I+3kA-k^2A^2+\mathrm{e}^{-kA}(4I+kA)]F(c_n,t_n+k)\,,
\end{split}
\end{equation}
where 
\begin{align*}
    a_n&= \mathrm{e}^{-\frac{k}{2}A}U_n - A^{-1}(\mathrm{e}^{-\frac{k}{2}A}-I)F(U_n,t_n)\\
    b_n&= \mathrm{e}^{-\frac{k}{2}A}U_n - A^{-1}(\mathrm{e}^{-\frac{k}{2}A}-I)F(a_n,t_n+\frac{k}{2})\\
    c_n&= \mathrm{e}^{-\frac{k}{2}A}a_n - A^{-1}(\mathrm{e}^{-\frac{k}{2}A}-I)\left[2F(b_n,t_n+\frac{k}{2})-F(U_n,t_n)\right]\,.  
\end{align*}
For convenience, we will utilize the following more compact form of the scheme, 
\begin{equation}
\begin{split}
\label{semidiscretecomp1}
    a_n&= \mathrm{e}^{-\frac{k}{2}A}U_n +\tilde{P}(kA)F(U_n,t_n)\\
    b_n&= \mathrm{e}^{-\frac{k}{2}A}U_n +\tilde{P}(kA)F(a_n,t_n+\frac{k}{2})\\
    c_n&= \mathrm{e}^{-\frac{k}{2}A}a_n \tilde{P}(kA)[2F(b_n,t_n+\frac{k}{2})-F(U_n,t_n)]\\
    U_{n+1} &= \mathrm{e}^{-kA}U_n +P_1(kA) F(U_n,t_n) +2P_2(kA)\left(F(a_n,t_n+\frac{k}{2})+F(b_n,t_n+\frac{k}{2})\right)\\
    &+ P_3(kA)F(c_n,t_n+k)\,, 
\end{split}
\end{equation}
with \begin{align}
\label{semidiscretecomp2}
    P_1(kA)&=  \frac{1}{k^2}(-A)^{-3}[-4I+kA+\mathrm{e}^{-kA}(4I+3kA+k^2A^2)]\nonumber \\
    P_2(kA)&=  \frac{1}{k^2}(-A)^{-3}[2I-kA-\mathrm{e}^{-kA}(2I+kA)]\\
    P_3(kA)&= \frac{1}{k^2}(-A)^{-3}[-4I+3kA-k^2A^2+\mathrm{e}^{-kA}(4I+kA)] \nonumber \\
    \tilde{P}(kA)&=- A^{-1}(\mathrm{e}^{-\frac{k}{2}A}-I).\nonumber
\end{align}

 Consider the decomposition $A = A_1 + A_2$, where $A_1$ and $A_2$ commute, we can rewrite \eqref{ADReqmatrix} in terms of $A_1$ and $A_2$ as follows.  Consider a new variable $V(t) = \mathrm{e}^{Bt}U(t)$ where $B$ is also a matrix, then $V^{'}(t) = \mathrm{e}^{Bt}U^{'}(t) + B\mathrm{e}^{Bt}U(t)$. Replacing $U^{'}(t)$ with the expression $-AU(t) + F(U(t),t)$ from \eqref{ADReqmatrix} we have $V^{'}(t) = \mathrm{e}^{Bt}F(U(t),t) - \mathrm{e}^{Bt}AU(t) + B\mathrm{e}^{Bt}U(t)$. Now, suppose $B$ commutes with $A$ then $\mathrm{e}^{Bt}AU(t) = A\mathrm{e}^{Bt}U(t)$ and $V^{'}(t) = \mathrm{e}^{Bt}F(U(t),t) + (B-A)\mathrm{e}^{Bt}U(t)$. Setting $B=A_1$ and using the decomposition of $A= A_1 + A_2$ we obtain the following semi-linear problem in the transformed variable $V$
 \begin{align}
     \label{vproblem}
     V^{'}(t) + A_2V &= \tilde{F}(V(t),t) \\
     V(0) &= U(0) \nonumber
 \end{align}
 with $ \tilde{F}(V,t)=\mathrm{e}^{A_1t}F(\mathrm{e}^{-A_1t}V,t)$. The solution variable $U(t)$ can then be obtained from the transformed variable $V(t)$ via the relation $U(t) = \mathrm{e}^{-A_1t}V(t)$. This reformulated semi-linear problem is the basis for deriving the dimensional splitting version of the class of ETDRK4 schemes. We accomplish this by applying the existing ETDRK4 scheme in \eqref{semidiscrete} to our reformulated problem \eqref{vproblem} and the formula $U(t) = \mathrm{e}^{-A_1t}V(t)$ to obtain a numerical solution for $U(t)$. We used a similar technique in \cite{asante2020second,asante2016dimensional} to develop the class of second-order dimensionally split ETD Runge-Kutta schemes. 
 
 Applying the ETDRK4 scheme to \eqref{vproblem} we obtain
 \begin{align}
 \label{vsemidiscrete}
    a_n^{*}&= \mathrm{e}^{-\frac{k}{2}A_2}V_n +\tilde{P}(kA_2)\tilde{F}(V_n,t_n) \nonumber\\
    b_n^{*}&= \mathrm{e}^{-\frac{k}{2}A_2}V_n +\tilde{P}(kA_2)\tilde{F}(a_n^{*},t_n+\frac{k}{2}) \nonumber\\
    c_n^{*}&= \mathrm{e}^{-\frac{k}{2}A_2}a_n^{*}+\tilde{P}(kA_2)\left[2\tilde{F}(b_n^{*},t_n+\frac{k}{2})-\tilde{F}(V_n,t_n)\right]\,\nonumber\\
    V_{n+1} &= \mathrm{e}^{-kA_2}V_n +P_1(kA_2) \tilde{F}(V_n,t_n) +  2P_2(kA_2)\left(\tilde{F}(a_n^{*},t_n+\frac{k}{2})+\tilde{F}(b_n^{*},t_n+\frac{k}{2})\right) \\
    &+ P_3(kA_2)\tilde{F}(c_n^{*},t_n+k). \nonumber
\end{align}

Now, using the fact that $V_n = \mathrm{e}^{A_1t_n}U_n$ and $\tilde{F}(V_n,t_n)=\mathrm{e}^{A_1t_n}F(\mathrm{e}^{-A_1t_n}V_n,t_n)$ and setting ${a_n = \mathrm{e}^{-A_1t_{n+\frac{1}{2}}}a_n^{*}}$, $b_n = \mathrm{e}^{-A_1t_{n+\frac{1}{2}}}b_n^{*}$, $c_n = \mathrm{e}^{-A_1t_{n+1}}c_n^{*}$, we can rewrite \eqref{vsemidiscrete} in terms of $U_n$,
\begin{equation}
\begin{split}
\label{vsemidiscrete2}
    a_n&= \mathrm{e}^{-\frac{k}{2}A_2}\mathrm{e}^{-\frac{k}{2}A_1}U_n +\tilde{P}(kA_2)\mathrm{e}^{-\frac{k}{2}A_1}F(U_n,t_n)\\
    b_n&= \mathrm{e}^{-\frac{k}{2}A_2}\mathrm{e}^{-\frac{k}{2}A_1}U_n +\tilde{P}(kA_2)F(a_n,t_n+\frac{k}{2})\\
    c_n&= \mathrm{e}^{-\frac{k}{2}A_2}\mathrm{e}^{-\frac{k}{2}A_1}a_n +\tilde{P}(kA_2)\left[2\mathrm{e}^{-\frac{k}{2}A_1}F(b_n,t_n+\frac{k}{2})-\mathrm{e}^{-kA_1}F(U_n,t_n)\right]\, \\
    U_{n+1} &= \mathrm{e}^{-kA_1}\mathrm{e}^{-kA_2}U_n +P_1(kA_2) \mathrm{e}^{-kA_1}F(U_n,t_n)  +2P_2(kA_2)\mathrm{e}^{-\frac{k}{2}A_1} G(a_n,b_n,t_n+\frac{k}{2})\\
    &+ P_3(kA_2)F(c_n,t_n+k) 
\end{split}
\end{equation}
with
\[G(a_n,b_n,t_n+\frac{k}{2})=F(a_n,t_n+\frac{k}{2})+F(b_n,t_n+\frac{k}{2}).\]
Equation \eqref{vsemidiscrete2} is the dimensional splitting formulation of the fourth-order ETD Runge-Kutta scheme first proposed by Cox and Matthews (\cite{Cox2002}, p. 436) for diagonal systems and extended to non-diagonal systems by Kassam and Trefethen (\cite{kassam2005fourth}, p. 1218).

\subsubsection{Pad\'e approximation of the matrix exponential }
Whereas the scheme in \eqref{vsemidiscrete2} is fully discretized and can be applied to solve the semi-linear problem in \eqref{ADReqmatrix}, it involves high powers of matrix inverses which can be expensive to compute, at best, or fail to exist, at worst, particularly in the cases of a coefficient matrix obtained from a PDE with Neumann boundary conditions. The scheme also requires an efficient implementation of the matrix exponential for reasonable computational time. To resolve this, rational functions are commonly used to approximate the matrix exponential and simplify the scheme in \eqref{vsemidiscrete2}. Pad\'e rational functions have been used successfully to develop the class of ETD-Pad\'e schemes leading to efficient second-order \cite{Kleefeld2012,Yousuf2012} and fourth-order ETD schemes \cite{Wade2009,yousuf2009efficient}.  Rational functions with real and distinct poles have also been used to develop the second-order schemes \cite{asante2020second,asante2016real}. In the discussion that follows we will employ the Pad\'e (2,2) rational function to approximate the matrix exponential in \eqref{vsemidiscrete2} and simplify the scheme. This same rational function was used to develop the fourth-order ETD-Pad\'e scheme in \cite{Wade2009} without dimensional splitting, referred to here as ETDRK4P22.  We denote the Pad\'e (2,2) approximation of $\mathrm{e}^{-kA}$ and $\mathrm{e}^{-\frac{k}{2}A}$ by $R_{2,2}(kA)$ and $\tilde{R}_{2,2}(kA)$, respectively, thus obtaining the representations
\begin{align}
    \mathrm{e}^{-kA}&\approx R_{2,2}(kA) =(12I-6kA+k^2A^2)(12I+6kA+k^2A^2)^{-1}\,, \label{pade1} \\
    \mathrm{e}^{-\frac{k}{2}A}&\approx \tilde{R}_{2,2}(kA) =(48I-12kA+k^2A^2)(48I+12kA+k^2A^2)^{-1}\label{pade2}\,.
\end{align}
Replacing the exponential functions in \eqref{vsemidiscrete2} by the corresponding rational approximations above, the functions $P_1, P_2, P_3,\tilde{P}$ can be expressed as the rational functions
\begin{align*}
    P_1(kA)&=  \frac{1}{k^2}(-A)^{-3}[-4I+kA+ R_{2,2}(kA_1)(4I+3kA+k^2A^2)]\nonumber \\
    &= k(2I-kA_2)(12I+6kA_2+k^2A^2)^{-1}\nonumber\\
    P_2(kA)&=  \frac{1}{k^2}(-A)^{-3}[2I-kA- R_{2,2}(kA_1)(2I+kA)]\\
    &=2k(12I + 6kA_2 + k^2A_2^2)^{-1}\\
    P_3(kA)&= \frac{1}{k^2}(-A)^{-3}[-4I+3kA-k^2A^2+R_{2,2}(kA_1)(4I+kA)]\\
    &= k(2I+kA_2)(12I+6kA_2+k^2A_2^2)^{-1} \nonumber\\
    \tilde{P}(kA_2) &= - A_2^{-1}(\tilde{R}_{2,2}(kA_2)-I)= 24k(48I+12kA+k^2A^2)^{-1}
\end{align*}
leading to the simplified scheme,
\begin{align}
    a_n&= \tilde{R}_{2,2}(kA_2)\tilde{R}_{2,2}(kA_1)U_n  +\tilde{P}(kA_2) \tilde{R}_{2,2}(kA_1)F(U_n,t_n) \label{ansplitsemi}\\
    b_n&= \tilde{R}_{2,2}(kA_2)\tilde{R}_{2,2}(kA_1)U_n  +\tilde{P}(kA_2) F(a_n,t_n+\frac{k}{2}) \label{bnsplitsemi}\\
    c_n&= \tilde{R}_{2,2}(kA_2)\tilde{R}_{2,2}(kA_1)a_n \nonumber\\
    &+\tilde{P}(kA_2) [2 \tilde{R}_{2,2}(kA_1)F(b_n,t_n+\frac{k}{2})- R_{2,2}(kA_1)F(U_n,t_n)] \label{cnsplitsemi} \\ 
    U_{n+1} &= R_{2,2}(kA_1)R_{2,2}(kA_2)U_n +P_1(kA_2) R_{2,2}(kA_1)F(U_n,t_n) \label{vsemidiscrete3} \\
    &+ 2P_2(kA_2)\tilde{R}_{2,2}(kA_1)G(a_n,b_n,t_n+\frac{k}{2})+P_3(kA_2)F(c_n,t_n+k) \nonumber
\end{align}
 We refer to the simplified scheme in \eqref{ansplitsemi}--\eqref{vsemidiscrete3} as the ETDRK4P22-IF scheme. 
 
 The use of Pad\'e (2,2) rational functions, which are well known to be A-acceptable, suggests that our scheme will have limited damping properties for non-smooth problems. For such problems, a few pre-smoothing steps involving a lower order scheme with an L-acceptable rational approximation to the exponential was recommended in \cite{Wade2009}. We investigate presmoothing the ETDRK4P22-IF scheme with the third order ETDRK scheme proposed in \cite{Wade2009}, which approximates the matrix exponentials in \eqref{semidiscretecomp1}--\eqref{semidiscretecomp2} with Pad\'e (0,3) rational functions. Note that this smoothing scheme is not implemented with dimensional splitting.

\subsubsection{Partial fraction decomposition of rational functions}

In implementing the ETDRK4P22-IF scheme, we utilize the partial fraction decomposition of the rational functions. This facilitates easy parallelization on multiple processors. The respective decomposition, whose complex poles and weights are similar to ones used in \cite{Yousuf2012}, are as follows:
\begin{equation}
\label{pfdR22}
R_{2,2}(z) = 1 + 2\mathrm{Re}\left(\frac{w_{11}}{z-c_1} \right), \ c_1 =-3+\mathrm{i}\sqrt{3},\ w_{11} = -\left(6 +\mathrm{i}\frac{18}{\sqrt{3}}\right)\,.
\end{equation}
\begin{equation}
\label{pfdRt22}
\tilde{R}_{2,2} = 1 + 4\mathrm{Re}\left(\frac{w_{11}}{z-c_2} \right), \ c_2 =-6+2\mathrm{i}\sqrt{3},\ w_{11} = -\left(6 +\mathrm{i}\frac{18}{\sqrt{3}}\right)\,.
\end{equation}
\begin{equation}
\label{pfdP1}
P_1(z) = 2k\mathrm{Re}\left(\frac{w_{21}}{z-c_1} \right), \ c_1 =-3+\sqrt{3}\mathrm{i},\ w_{21} = -\left(\frac{1}{2}+\mathrm{i}\frac{5\sqrt{3}}{6}\right)\,.
\end{equation}
\begin{equation}
\label{pfdP2}
P_2(z) = 4k\mathrm{Re}\left(\frac{w_{31}}{z-c_1} \right), \ c_1 =-3+\sqrt{3}\mathrm{i},\ w_{31} = -\mathrm{i}\frac{\sqrt{3}}{6}.
\end{equation}
\begin{equation}
\label{pfdP3}
P_3(z) = 2k\mathrm{Re}\left(\frac{w_{41}}{z-c_1} \right), \ c_1 =-3+\mathrm{i}\sqrt{3},\ w_{41} = \frac{1}{2}+\mathrm{i}\frac{\sqrt{3}}{6}.
\end{equation}
\begin{equation}
\label{pfdPtilde}
\tilde{P}(z) = 48k\mathrm{Re}\left(\frac{w_{51}}{z-c_2} \right), \ c_2 =-6+2\mathrm{i}\sqrt{3},\ w_{51} = -\mathrm{i}\frac{\sqrt{3}}{12}.
\end{equation}

Partial fraction decomposition of rational approximations to the exponential was also used to develop the ETD Crank Nicolson scheme \cite{Kleefeld2012} and ETDRDP scheme \cite{asante2016real}, where in the later all poles and weights are real. 

% The implementation of the new ETDRK4P22-IF scheme requires further simplification of \eqref{vsemidiscrete3} using the partial fraction decomposition of the Pad\'e rational functions derived above. It is instructive to first describe the implementation of the ETDRK4P22 scheme as some of the expressions will be useful in shortening the discussion of the implementation of the dimensional splitting scheme. The simplified ETDRK4P22 scheme obtained at end of the next section is identical to the scheme developed in \cite{Wade2009}.
\subsubsection{Implementation of ETDRK4P22 scheme}
Before laying out the implementation steps for the new ETDRK4P22-IF scheme, it is instructive to first describe the implementation of the unsplit ETDRK4P22 scheme as some expressions will be carried over to shorten the discussion.

 The ETDRK4P22 scheme expressed in terms of the respective rational approximations of the matrix exponential is
\begin{align}
\label{etdrk422unsplit}
    a_n&= \tilde{R}_{2,2}(kA)U_n  +\tilde{P}(kA)F(U_n,t_n) \nonumber \\
    b_n&= \tilde{R}_{2,2}(kA)U_n  +\tilde{P}(kA)F(a_n,t_n+\frac{k}{2})\nonumber \\
    c_n&= \tilde{R}_{2,2}(kA)a_n+\tilde{P}(kA) [2F(b_n,t_n+\frac{k}{2})-F(U_n,t_n)]\, \nonumber\\
    U_{n+1} &= R_{2,2}(kA)U_n + P_1(kA)F(U_n,t_n) + 2P_2(kA)G(a_n,b_n, t_n + \frac{k}{2})\\
    &+ P_3(kA)F(c_n,t_n+k)\,, \nonumber
\end{align} To implement this scheme, we first need to compute $a_n$, $b_n$ and $c_n$. Substituting the decomposition in \eqref{pfdRt22} and \eqref{pfdPtilde} into \eqref{etdrk422unsplit}, $a_n$ becomes
\begin{align}
    a_n&= \tilde{R}_{2,2}(kA)U_n + \tilde{P}(kA) F(U_n,t_n) \nonumber\\
       &= \left( I + 4\mathrm{Re}\left(w_{11}(kA - c_2I)^{-1}\right) \right) U_n+ 48k\mathrm{Re}\left( w_{51}(kA-c_2I)^{-1} \right) F(U_n,t_n)\nonumber\\
       &=  U_n + 2\mathrm{Re}\left[ (kA - c_2I)^{-1}( 2w_{11} U_n )\right]+2k\mathrm{Re}\left[(kA - c_2I)^{-1}(24w_{51}F(U_n,t_n)) \right] \nonumber \\
       &= U_n +2\mathrm{Re}(a_{n1})\,,  \label{etdrk4p22an}
\end{align}
where $ (kA - c_2I)a_{n1}= 2w_{11} U_n+ k24w_{51}F(U_n,t_n)$. For $b_n$ we get
\begin{align}
    b_n&= \tilde{R}_{2,2}(kA)U_n +   \tilde{P}(kA) F(a_n,t_n+\frac{k}{2})\nonumber \\
       &=  U_n + \mathrm{Re}\left[ (kA - c_2I)^{-1}( 4w_{11} U_n )\right]+k\mathrm{Re}\left[(kA -c_2I)^{-1}(48w_{51}F(a_n,t_n+\frac{k}{2})) \right] \nonumber \\
       &= U_n + 2\mathrm{Re}(b_{n1})\,, \label{etdrk4p22bn}
\end{align}
where $(kA - c_2)b_{n1}=2w_{11} U_n + 24w_{51}kF(a_n,t_n+\frac{k}{2})$. Finally for $c_n$ we have, 
\begin{align}
    c_n&= \tilde{R}_{2,2}(kA)a_n +   \tilde{P}(kA) [2F(b_n,t_n+\frac{k}{2})-F(U_n,t_n)] \nonumber\\
       &= \left( I + 4\mathrm{Re}\left(w_{11}(kA - c_2I)^{-1}\right) \right) a_n+ 48k\mathrm{Re} \left( w_{51}(kA-c_2I)^{-1} \right)\times \nonumber\\
       &[2F(b_n,t_n+\frac{k}{2})-F(U_n,t_n)] \nonumber \\
       &=  a_n + 2\mathrm{Re}\left[ (kA - c_2I)^{-1}\left( 2w_{11} a_n +24w_{51}k[2F(b_n,t_n+\frac{k}{2})-F(U_n,t_n)]\right) \right] \nonumber \\
       &= a_n + 2\mathrm{Re}(c_{n1})\,, \label{etdrk4p22cn}
\end{align}
where $ (kA - c_2I)c_{n1}= 2w_{11} a_n + 24w_{51}k(2F(b_n,t_n+\frac{k}{2})-F(U_n,t_n))$.

For the main scheme we make use of the decomposition \eqref{pfdR22}, \eqref{pfdP1}--\eqref{pfdP3} to obtain,
\begin{align}
    U_{n+1} &=  R_{2,2}(kA)U_n + P_1(kA)F(U_n,t_n) + 8k\mathrm{Re}( w_{31}(kA-c_1I)^{-1})G(a_n,b_n, t_n +\frac{k}{2}) \nonumber \\
    &+ 2k\mathrm{Re}(w_{41}(kA-c_1I)^{-1})F(c_n,t_n+k)\nonumber \\
        &= \left( I + 2\mathrm{Re}\left(w_{11}(kA - c_1I)^{-1}\right) \right) U_n+ 2k\mathrm{Re} \left( w_{21}(kA-c_1I)^{-1} \right) F(U_n,t_n) \nonumber \\
        &+ 2k \mathrm{Re} \left[ (kA-c_1I)^{-1} (4w_{31}G(a_n,b_n,t_n+\frac{k}{2})+w_{41}F(c_n,t_{n}+k))\right]\nonumber \\
        &= U_n + 2\mathrm{Re}(U_{n1}) \,,\label{etdrk4p22dn}
\end{align}
where $ (kA-c_1I)U_{n1}=w_{11}U_n + kw_{21} F(U_n,t_n) +4w_{31}kG(a_n,b_n,t_n+\frac{k}{2})+w_{41}kF(c_n,t_{n}+k))$.
Putting all the results together, the ETDRK4P22 scheme can be implemented in serial with the following steps:

\begin{enumerate}
    \item Solve: $(kA-c_2I)a_{n1}= 2w_{11}u_n + 24w_{51}kF(u_n,t_n)$
    \item Set: $a_n = u_n + 2\mathrm{Re}(a_{n1})$
    \item Solve: $(kA-c_2I)b_{n1} =2w_{11}u_n + 24w_{51}kF(a_n,t_n+\frac{k}{2})$
    \item Set: $b_n = u_n + 2\mathrm{Re}(b_{n1})$
    \item Solve:  $(kA-c_2I)c_{n1} = 2w_{11}a_n+ 24w_{51}k[2F(b_n,t_n+\frac{k}{2})- F(u_n,t_n)] $
    \item Set:  $c_n = a_n +2\mathrm{Re}(c_{n1})$
    \item Solve: $(kA-c_1I)u_{n1} = w_{11}u_n + w_{21}kF(u_n,t_n)+ 4w_{31}kG(a_n,b_n,t_n+\frac{k}{2}) + w_{41}kF(c_n,t_n+k)$.
    \item Set: $u_{n+1}= u_n + 2\mathrm{Re}(u_{n1})$.
\end{enumerate}

\subsubsection{Implementation of ETDRK4P22-IF scheme}
\label{finalscheme}
Following from the derivation above, we lay out the implementation of the dimensionally split scheme begining from the intermediate solutions $a_n$, $b_n$ and $c_n$. For $a_n$ we make use of \eqref{etdrk4p22an} to obtain,
\begin{align*}
    a_n &= \tilde{R}_{2,2}(kA_1)\left[ \tilde{R}_{2,2}(kA_2)U_n + \tilde{P}(kA_2)F(U_n,t_n)\right]\\
        &= \tilde{R}_{2,2}(kA_1)( u_n +2\mathrm{Re}(a_{n1}))
        \intertext{where  $ (kA_2 - c_2I)a_{n1}= 2w_{11} U_n+ k24w_{51}F(U_n,t_n)$. Setting $a_{n2}= U_n +2\mathrm{Re}(a_{n1})$ and applying the partial fraction decomposition of  $\tilde{R}_{2,2}(kA_1)$, we have}
   a_n &= \left( I + 4\mathrm{Re}\left(w_{11}(kA_1 - c_2I)^{-1}\right) \right) a_{n2}\\
        &= a_{n2} + 2\mathrm{Re}((kA_1 - c_2I)^{-1}(2w_{11}a_{n2}))\\
        &= a_{n2} + 2\mathrm{Re}(a_{n3}),
\end{align*}
with $(kA_1 - c_2I)a_{n3}= 2w_{11}a_{n2}$. Similarly for $b_n$  we use \eqref{etdrk4p22bn} to obtain
\begin{align*}
    b_n &=  \tilde{R}_{2,2}(kA_1)\tilde{R}_{2,2}(kA_2)U_n  +\tilde{P}(kA_2) F(a_n,t_n+\frac{k}{2})\\
        &= \tilde{R}_{2,2}(kA_1)( U_n +2\mathrm{Re}(b_{n1})) + 2\mathrm{Re}(b_{n2})),
        \intertext{where  $ (kA_2 - c_2I)b_{n1}= 2w_{11} U_n$ and $ (kA_2 - c_2I)b_{n2}= 24kw_{51} F(a_n,t_n +\frac{k}{2})$. Setting $b_{n3}= U_n +2\mathrm{Re}(b_{n1})$ and applying the partial fraction decomposition of  $\tilde{R}_{2,2}(kA_1)$ we see that}
    b_n &= \left( I + 4\mathrm{Re}\left(w_{11}(kA_1 - c_2I)^{-1}\right) \right) b_{n3}  + 2\mathrm{Re}(b_{n2})\\
        &= b_{n3} + 2\mathrm{Re}((kA_1 - c_2I)^{-1}(2w_{11}b_{n3})) + 2\mathrm{Re}(b_{n2})\\
        &= b_{n3} + 2\mathrm{Re}(b_{n4}) +  2\mathrm{Re}(b_{n2}),
\end{align*}
with $(kA_1 - c_2I)b_{n4}= 2w_{11}b_{n3}$. Similarly from \eqref{etdrk4p22cn},
\begin{align*}
    c_n &= \tilde{R}_{2,2}(kA_1)[\tilde{R}_{2,2}(kA_2)a_n+2\tilde{P}(kA_2)F(b_n,t_n+\frac{k}{2})]- R_{2,2}(kA_1)\tilde{P}(kA_2)F(U_n,t_n) ]\\
        &= \tilde{R}_{2,2}(kA_1)c_{n1}^{*}- R_{2,2}(kA_1)c_{n2}^{*},
\end{align*} 
where $c_{n1}^{*}=\tilde{R}_{2,2}(kA_2)a_n+2\tilde{P}(kA_2)F(b_n,t_n+\frac{k}{2})$ and $c_{n2}^{*}=\tilde{P}(kA_2)F(U_n,t_n)$. Using the respective partial fraction decomposition we get,
    $$c_n = c_{n1}^{*} + 2\mathrm{Re}(c_{n3}) - [c_{n2}^{*} + 2\mathrm{Re}(c_{n4})]$$ where  $ (kA_1 - c_2I)c_{n3}= 2w_{11} c_{n1}^{*}$ and  $ (kA_1 - c_1I)c_{n4}= w_{11} c_{n2}^{*}$. Now it remains to find $c_{n1}^{*}$ and $c_{n2}^{*}$. We know that
    $$c_{n1}^{*} = \tilde{R}_{2,2}(kA_2)a_n+2\tilde{P}(kA_2)F(b_n,t_n+\frac{k}{2})=a_n + 2\mathrm{Re}(c_{n1})$$ with $(kA_2-c_2I)c_{n1}= 2w_{11}a_n + 48w_{51}kF(b_n,t_n+\frac{k}{2})$. Similarly we obtain, $$c_{n2}^{*} =\tilde{P}(kA_2)F(U_n,t_n)= 2\mathrm{Re}(c_{n2})\,,$$ 
 where $(kA_2-c_2I)c_{n2}=24w_{51}kF(U_n,t_n)$. Now for the main scheme, set $U_{n1}^{*} = R_{2,2}(kA_2)U_n +P_1(kA_2) F(U_n,t_n)$, $U_{n2}^{*}=  2P_2(kA_2)G(a_n,b_n,t_n+\frac{k}{2})$ and $U_{n3}^{*}=P_3(kA_2) F(c_n,t_n+k)$, then we can express $U_{n+1}$ as 
{\footnotesize \begin{align*}
    U_{n+1} &= R_{2,2}(kA_1)U_{n1}^{*} + \tilde{R}_{2,2}(kA_1)U_{n2}^{*} + U_{n3}^{*}\\
           &=  \left( I + 2\mathrm{Re}\left(w_{11}(kA_1 - c_1I)^{-1}\right) \right)U_{n1}^{*} +  \left( I + 4\mathrm{Re}\left(w_{11}(kA_1 - c_2I)^{-1}\right) \right)U_{n2}^{*} +  U_{n3}^{*} \\
           &= U_{n1}^{*} + U_{n2}^{*} + 2\mathrm{Re}( (kA_1-c_1I)^{-1}(w_{11}U_{n1}^{*}) ) + 2\mathrm{Re}((kA_1-c_2I)^{-1}(2w_{11}U_{n2}^{*})) +  U_{n3}^{*} \\
           &= U_{n1}^{*} + U_{n2}^{*} + 2\mathrm{Re}( U_{n4} ) + 2\mathrm{Re}(U_{n5}) + U_{n3}^{*} \,,
\end{align*}}
where $(kA_1-c_1I)U_{n4}= w_{11}U_{n1}^{*}$ and $(kA_1-c_2I)U_{n5}=2w_{11}U_{n2}^{*}$. Now from \eqref{etdrk4p22dn} $U_{n1}^{*} = U_n + 2\mathrm{Re}(U_{n1})$ where $(kA_2-c_1I)U_{n1}=w_{11}U_n + kw_{21} F(U_n,t_n).$  Similarly,  $U_{n2}^{*}= 2\mathrm{Re}(U_{n2})$ where $(kA_2-c_1I)U_{n2}= 4w_{31}kG(a_n,b_n,t_n+\frac{k}{2}).$ Finally, $U_{n3}^{*}=2\mathrm{Re}(U_{n3})$ where $(kA_2-c_1I)U_{n3}=w_{41}kF(c_n,t_n+k)$.

\subsubsection*{Algorithm for 2D ETDRK4P22-IF}
The full ETDRK4P22-IF scheme can be implemented in parallel with up to three processors in the following steps. 
\begin{enumerate}
    \item Solve: $(kA_2-c_2I)a_{n1}= 2w_{11}U_n+24w_{51}kF(U_n,t_n)$
    \item Set:  $a_{n2}= U_n +2\mathrm{Re}(a_{n1})$
    \item Solve: $(kA_1-c_2I)a_{n3}= 2w_{11}a_{n2}$
    \item Set: $a_n = a_{n2}+2\mathrm{Re}(a_{n3})$.
    \item Solve Processor 1: $(kA_2-c_2I)b_{n1}= 2w_{11}U_n$
    \item Solve Processor 2: $(kA_2-c_2I)b_{n2}= 24w_{51}kF(a_n,t_n+\frac{k}{2})$
    \item Set:  $b_{n3}= U_n +2\mathrm{Re}(b_{n1})$
    \item Solve: $(kA_1-c_2I)b_{n4}= 2w_{11}b_{n3}$
    \item Set: $b_n = b_{n3}+2\mathrm{Re}(b_{n4})+ 2\mathrm{Re}(b_{n2})$.
    \item Solve Processor 1: $(kA_2-c_2I)c_{n1}= 2w_{11} a_n+ 48w_{51}kF(b_n,t_n+\frac{k}{2})$
    \item Solve Processor 2: $(kA_2-c_2I)c_{n2}= 24w_{51}kF(U_n,t_n)$
    \item Set: $c_{n1}^{*}=a_n+2\mathrm{Re}(c_{n1})$ and $c_{n2}^{*}=2\mathrm{Re}(c_{n2})$
    \item Solve Processor 1: $(kA_1-c_2I)c_{n3}= 2w_{11}c_{n1}^{*}$
    \item Solve Processor 2: $(kA_1-c_1I)c_{n4}= w_{11}c_{n2}^{*}$
    \item Set: $c_{n} = c_{n1}^{*} + 2\mathrm{Re}(c_{n3}) - [c_{n2}^{*} + 2\mathrm{Re}(c_{n4})]$
    \item Solve Processor 1: $(kA_2-c_1I)U_{n1} =  w_{11}U_n+ w_{21}kF(U_n,t_n) $
    \item Solve Processor 2: $(kA_2-c_1I)U_{n2} = 4w_{31}kG(a_n,b_n,t_n+\frac{k}{2})$
    \item Solve Processor 3: $(kA_2-c_1I)U_{n3} = w_{41}kF(c_n,t_n+k)$ 
    \item Set: $u_{n1}^{*} = u_n + 2 \mathrm{Re}(U_{n1})$, $U_{n2}^{*} = 2\mathrm{Re}(U_{n2})$  $U_{n3}^{*} = 2\mathrm{Re}(U_{n3})$
    \item Solve Processor 1: $(kA_1-c_1I)U_{n4} = w_{11}U_{n1}^{*}$  
    \item Solve Processor 2: $(kA_1-c_2I)U_{n5} = 2w_{11}U_{n2}^{*}$ 
    \item Set: $U_{n+1}= U_{n1}^{*} + U_{n2}^{*}+ U_{n3}^{*} + 2\mathrm{Re}(U_{n4}) + 2\mathrm{Re}(U_{n5})$.
\end{enumerate}

For non-smooth problems the following third order scheme, developed in \cite{Wade2009}, can be used for presmoothing. We show in our numerical experiments that applying a few steps of this presmoothing scheme followed by ETDRK4P22-IF preserves the fourth order accuracy of ETDRK4P22-IF. 
\begin{enumerate}
    \item Solve Processor 1: $(kA-f_1I)a_{n1}=2s_{11}u_n + ks_{51}F(u_n,t_n) $
    \item Solve Processor 2: $(kA-f_2I)a_{n2}=2s_{12}u_n + ks_{52}F(u_n,t_n)$
    \item Set $a_n = a_{n1}+2\mathrm{Re}(a_{n2})$.
    \item Solve Processor 1: $(kA-f_1I)b_{n1} = 2s_{11}u_n + ks_{51}F(a_n,t_n+\frac{k}{2})$
    \item Solve Processor 2: $(kA-f_2I)b_{n2}=  2s_{12}u_n + ks_{52}F(a_n,t_n+\frac{k}{2})$
    \item Set $b_n = b_{n1}+2\mathrm{Re}(b_{n2})$.
    \item Solve Processor 1:  $(kA-f_1I)c_{n1} = 2s_{11}a_n + ks_{51}(2F(b_n,t_n+\frac{k}{2})-F(u_n,t_n)) $
    \item Solve Processor 2: $(kA-f_2I)c_{n2}= 2s_{12}a_n + ks_{52}(2F(b_n,t_n+\frac{k}{2})-F(u_n,t_n))$.
    \item  Set $c_n = c_{n1} +2\mathrm{Re}(c_{n2})$.
    \item Solve Processor 1: $(kA-e_1I)u_{n1} = s_{11}u_n + ks_{21}F(u_n,t_n)+ 2ks_{31}G(a_n,b_n,t_n+\frac{k}{2}) + ks_{41}F(c_n,t_n+k)$. 
    \item Solve Processor 2 $(kA-e_2I)u_{n2}= s_{12}u_n + ks_{22}F(u_n,t_n)+ 2ks_{32}G(a_n,b_n,t_n+\frac{k}{2})) + ks_{42}F(c_n,t_n+k)$.
    \item Set $u_{n+1} = u_{n1} + 2\mathrm{Re}(u_{n2})$.
\end{enumerate}
The poles and weights are given by
\begin{align*}
    f_1 &= -3.19214327596664; \ f_2= -1.40392836201668 + 3.61467898890404\mathrm{i}\\
     e_1 &= 0.5f_1; \ e2= 0.5f_2; \ s_{11}=\frac{6}{|e_1-e_2|^2}; s_{12}=\frac{-3\mathrm{i}}{(\mathrm{Im}(e_2))(e_2-e_1)}\\
     s_{21}&=\frac{1-e_1}{|e_1-e_2|^2}; s_{22}=\frac{\mathrm{i}(e_2-1)}{2(\mathrm{Im}(e_2))(e_2-e_1)};\ s_{31}=\frac{1+e_1}{|e_1-e_2|^2}; s_{32}=\frac{-\mathrm{i}(e_2+1)}{2(\mathrm{Im}(e_2))(e_2-e_1)} \\
      s_{41}&=\frac{1+e_1^2}{|e_1-e_2|^2}; s_{42}=\frac{-\mathrm{i}(e_2^2+1)}{2(\mathrm{Im}(e_2))(e_2-e_1)};\ s_{51}=\frac{24+6f_1+f_1^2}{|f_1-f_2|^2}; s_{52}=\frac{-\mathrm{i}(24+6f_2+f_2^2)}{2(\mathrm{Im}( f_2))(f_2-f_1)}.
\end{align*}
%\subsection{Implementation and parallelization of the fully discrete scheme}
%

\section{Results and Discussion}
\label{sec:experiments}
We investigate the order of accuracy and efficiency of the proposed scheme ETDRK4P22-IF, which applies integrating factor-type dimensional splitting to improve the run time of the existing fourth-order ETDRK4P22 scheme developed in \cite{Wade2009}. We also compare the accuracy of the new scheme with the existing second-order ETDRDP-IF developed in \cite{asante2020second} with dimensional splitting as well as the fourth-order semi-implicit backward differentiation scheme (SBDF4) developed in \cite{ascher1995implicit}. SBDF4 is an IMEX scheme which employs a fourth-order backward differentiation formula to discretize the linear part of \eqref{ADReqmatrix} and a fourth order Adam's Bashford scheme to discretize the non-linear part. The resulting scheme is of the form
\begin{equation}
\begin{split}
    \label{sbdf4}
    (25I+13kA)U_{n+1} &= 48U_n-36U_{n-1}+16U_{n-2}-3U_{n-3}\\
    &+ k\left(48F(U_{n},t_n)-72kF(U_{n-1},t_{n-1})\right.\\
    &+\left.48kF(U_{n-2},t_{n-2})-12kF(U_{n-3},t_{n-3})\right).
\end{split}
\end{equation}
The initial stages of the scheme $U_1,U_2,U_3$ are computed using the first-order SBDF scheme \[(I+k_0A)U_{n+1}= U_n + k_0F(U_n,t_n)\] with a very small time step of $k_0 = \frac{k}{2000}$, as was previously suggested in \cite{ascher1995implicit,kassam2005fourth}.

We evaluate the performance of the schemes on problems in one and two spatial dimensions as well as for scalar and systems of reaction-diffusion equations with Dirichlet and Neumann boundary conditions. The error, $E(k)$ at time step $k$, is measured with the $L_{\infty}$ norm and the estimated rate of convergence is calculated with the formula
\begin{equation}
\label{eqn:conv1}
    p = \frac{\log[E(k)/E(k/2)]}{\log(2)}
\end{equation}
For problems with exact solution, $E(k)=\norm{U(k)-\hat{U}(k)}$, where $U(k),\hat{U}(k)$ represent the exact and approximate grid solutions at the final time $T$ using a temporal step size of $k$. All spatial derivatives are approximated with a fourth-order finite difference scheme on a mesh with size $h$ which enabled us to derive the correct order of convergence of the fully discrete scheme by setting $k=h$. For problems where the exact solution is not known, we only evaluated the convergence of the time discretization by fixing the spatial mesh size, $h$. The errors are computed by using the numerical solution on the next fine mesh as the reference solution. For a grid refinement study with temporal step sizes $k,k/2, k/4$, the formula
\begin{equation}
\label{eqn:conv2}
\tilde{p} = \frac{\log[\tilde{E}(k)/\tilde{E}(k/2)]}{\log(2)}
\end{equation}
is used,
where $\tilde{E}(k) = \norm{\hat{U}(k)-\hat{U}(k/2)}$ has been shown in \cite{leveque2007finite} to determine the correct rate of convergence of the scheme for sufficiently small $k$. Each numerical scheme evaluated was used to calculate its own reference solution. Numerical experiments were run in MATLAB R2020b on a THINKMATE Work Station with Intel (R) Core(TM) i9-10900k CPU @ 3.70 GHz with 128 GB RAM. All codes used to generate convergence tables and figures have been made publicly available at \url{https://github.com/easantea/ETDRK4P22-IF}.

\subsection{Model problem with Dirichlet boundary}
\label{Ex:Model1}
In order to investigate the order of convergence of the proposed scheme, we consider the following two-dimensional reaction-diffusion equation with exact solution
\begin{equation}
\frac{\partial u}{\partial t} = \Delta u - u\,,\quad  -\frac{\pi}{2}<x\,,\quad y<\frac{\pi}{2}\,,\quad  t\in [0, T]\,.
\end{equation}
with homogeneous Dirichlet boundary conditions and initial condition 
\[ u(x,y,0)=\cos(x)\cos(y). \] 
The exact solution is $u(x,y,t) = \mathrm{e}^{-3t}\cos(x)\cos(y)$.

\begin{table}[h]
 \begin{center}
  \caption{Convergence table showing fourth-order convergence of the new ETDRK4P22-IF compared with the existing ETDRK4P22 scheme for the model test problem with exact solution at final time $T=1$}
 \resizebox{\textwidth}{!}{
 \begin{tabular}{|c|c|c|c|c|c|c|c|}
 \hline
      &  & \multicolumn{3}{|c|}{ETDRK4P22-IF} & \multicolumn{3}{|c|}{ETDRK4P22} \\
  $k$ & $h$ & Error& $p$ & Cpu Time (sec) & Error & $p$ & Cpu Time (sec) \\

 \hline
 0.10 & 0.08 & 1.639$\times 10^{-7}$ &-- &0.06& 9.069$\times 10^{-7}$ & -- &0.02\\
 \hline
 $0.05$  & 0.04 & 1.0805 $\times 10^{-8}$ &3.92& 0.08 & 5.6131$\times 10^{-8}$ & 4.01 &0.35\\
\hline
 $0.025$  & 0.02 & 6.958$\times 10^{-10}$ &3.96& 0.65 &3.496$\times 10^{-9}$ & 4.01 &5.21\\
\hline
 $0.0125$  & 0.010 & 4.456$\times 10^{-11}$ &3.96&3.72 &2.1391$\times 10^{-10}$ & 4.03 &79.70\\
 \hline
  \end{tabular}
  }
  \end{center}
    \label{tab:convergenceEx1}
  \end{table}

\begin{figure}[ht]
   \centering
    \includegraphics[width=0.95\textwidth]{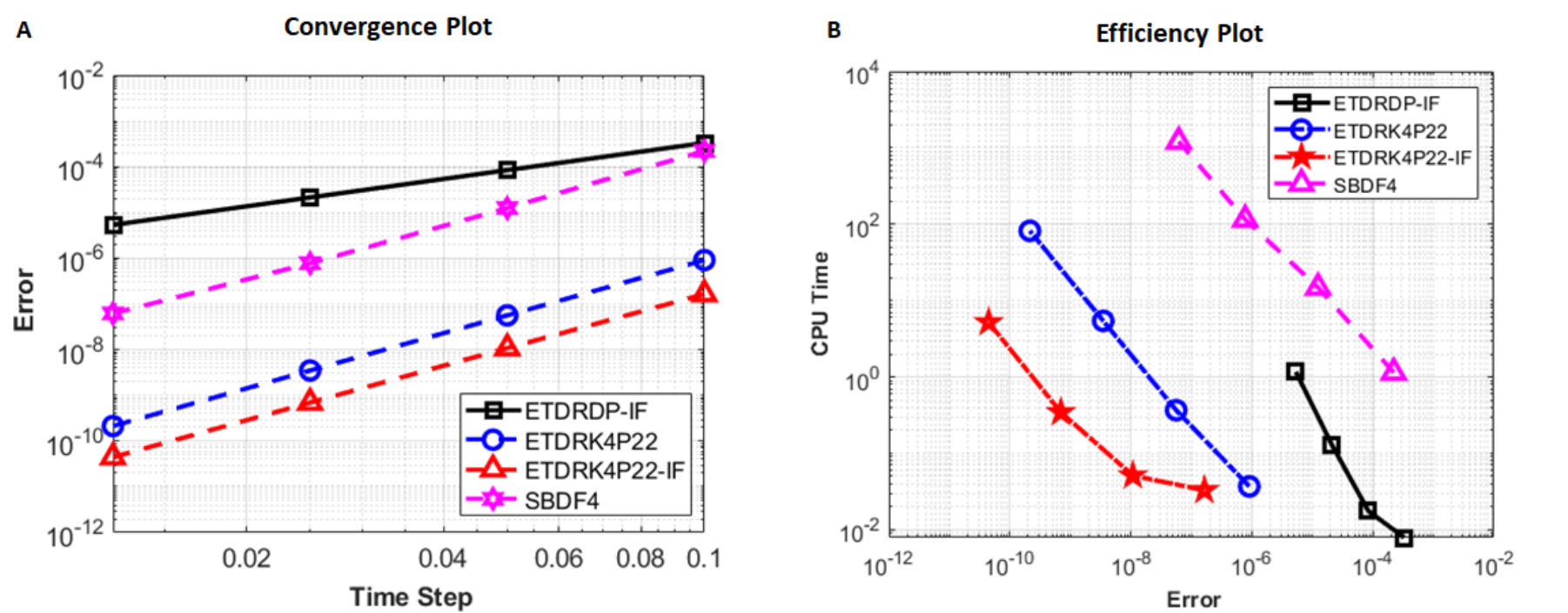}       
 \caption{\textbf{Convergence and efficiency plots for Example~\ref{Ex:Model1}}. A) Fourth-order convergence and B) Efficiency of the new fourth-order ETDRK4P22-IF scheme compared with the existing fourth-order ETDRK4P22, second-order ETDRDP-IF and the fourth order SBDF4 schemes.}
 \label{fig:efficiencyEx1}
 \end{figure}

We discretize the spatial derivative in each dimension using the following fourth-order central difference approximation described in Section \ref{sec:space}, where a fourth degree polynomial extrapolation is used to compute the derivatives near the boundary nodes ($j=0,j=m+1$). Each dimension of our spatial domain is descritized with $m+2$ nodes with $h=\frac{\pi}{m+1}$. Linear systems in the numerical scheme are solved using Gaussian elimination with LU-decomposition. From Table~\ref{tab:convergenceEx1} and Fig.~\ref{fig:efficiencyEx1}A we see that proposed dimensional splitting scheme ETDRK4P22-IF is fourth-order accurate. For a serial implementation of both schemes, ETDRK4P22-IF is 20 times faster than ETDRK4P22 at the finest spatial mesh size ($h=0.010$) and temporal step size ($k=0.0125$) (Table~\ref{tab:convergenceEx1}). Interestingly, the new ETDRK4P22-IF is more accurate than the existing ETDRK4P22 scheme across all temporal and spatial step sizes examined (Fig~\ref{fig:efficiencyEx1}B, Table~\ref{tab:convergenceEx1}). Whereas ETDRK4P22-IF requires more CPU time per time step compared to the second-order ETDRDP-IF, the errors for ETDRK4P22-IF are several orders of magnitude smaller than those for ETDRDP-IF (Fig.~\ref{fig:efficiencyEx1}B, Table~\ref{tab:AppEx1RDPSBDF4}).

For all time steps examined, the SBDF4 scheme is less accurate and requires substantially greater CPU time compared with ETDRK4P22-IF (Fig.~\ref{fig:efficiencyEx1}B, Table~\ref{tab:AppEx1RDPSBDF4}). Most of this time is spent computing the starting values of the multistep method. Increasing the time step used in computing the starting values decreased the CPU time but also worsened the accuracy of the scheme.

% \begin{enumerate}
%     \item Use Fourth-order centered finite difference discretization of spatial derivative.
%     \item Apply Forth-order ETDRK4Pade-SIF.
%     \item Compare dimensional splitting scheme with unsplit scheme.
%     \item Compare parallel implementation with serial. 

% \end{enumerate}

\subsection{Model problem with Neumann boundary}
\label{Ex:Model2}
Next, we investigate the order of convergence of the proposed scheme for a scalar reaction-diffusion problem with homogeneous Neumann boundary conditions. The model problem is given by
\begin{equation}
\frac{\partial u}{\partial t} = \Delta u -u\,,\quad -\pi<x\,,\quad y<\pi\,,\quad t\in [0, T]\,.
\end{equation}
with homogeneous Neumann boundary conditions and initial condition 
\[ u(x,y,0)=\cos(x)\cos(y)\,\textcolor{red}{.}\] 
The exact solution is $u(x,y,t) = \mathrm{e}^{-3t}\cos(x)\cos(y)$.

 \begin{table}[h]
 \begin{center}
      \caption{Convergence table showing fourth-order convergence of the new ETDRK4P22-IF scheme compared with the existing ETDRK4P22 scheme for the model test problem in Example~\ref{Ex:Model2} with $T=1$.}
 \resizebox{\textwidth}{!}{
 \begin{tabular}{|c|c|c|c|c|c|c|c|}
 \hline
      &  & \multicolumn{3}{|c|}{ETDRK4P22-IF} & \multicolumn{3}{|c|}{ETDRK4P22} \\
  $k$ & $h$ & Error& $p$ & Cpu Time (sec) & Error & $p$ & Cpu Time (sec) \\

 \hline
 0.10& 0.1 & 1.0836$\times 10^{-5}$ &-- &0.003& 1.1580$\times 10^{-5}$ & -- &0.002\\
 \hline
 $0.05$  & 0.05 & 6.8127 $\times 10^{-7}$ &3.99& 0.012 & 7.2661$\times 10^{-7}$ & 3.99 &0.042\\
\hline
 $0.025$  & 0.025 & 4.2638$\times 10^{-8}$ &4.00& 0.081 &4.5439$\times 10^{-8}$ & 4.00 &0.732\\
\hline
 $0.0125$  & 0.0124 & 2.6657$\times 10^{-9}$ &4.00&1.961 &2.8397$\times 10^{-9}$ & 4.00 & 10.703\\
 \hline
  \end{tabular}
  }
  \end{center}
    \label{tab:convergenceEx2}
  \end{table}

    \begin{figure}[h]
   \centering
    \includegraphics[width=0.95\textwidth]{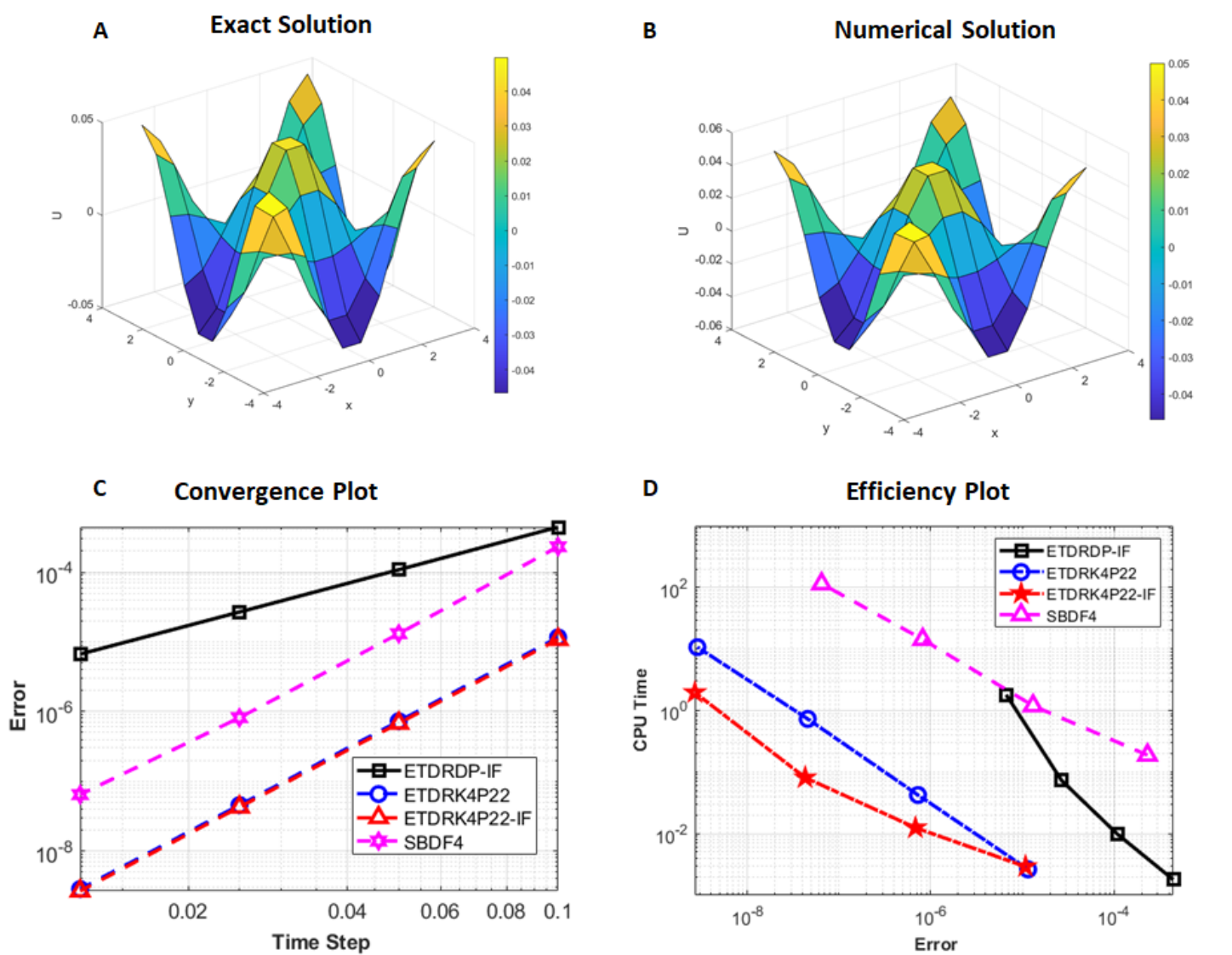}       
 \caption{\textbf{Convergence and efficiency plots for Example~\ref{Ex:Model2}}. A) Exact solution, B) Numerical solution and log-log plots showing C) Fourth-order convergence and D) Efficiency of the new fourth-order ETDRK4P22-IF scheme compared with the existing fourth-order ETDRK4P22, second-order ETDRDP-IF and fourth order SBDF4 schemes. }
 \label{fig:efficiencyEx2}
 \end{figure}

We partition the spatial domain into $m+2$ nodes with $h=\frac{2\pi}{m+1}$. The second derivative in each dimension is discretized using the fourth-order central difference approximation described in Section \ref{sec:space}. The resulting linear systems in the numerical scheme is solved using Gaussian elimination with LU-decomposition. The exact and numerical solution profiles are shown in Fig.~\ref{fig:efficiencyEx2}(A-B). The proposed scheme converges with fourth-order accuracy (Table~\ref{tab:convergenceEx2}, Fig~\ref{fig:efficiencyEx2}C). In a serial implementation, ETDRK4P22-IF is 5 times faster than ETDRK4P22 at the finest spatial mesh size ($h=0.0124$) and temporal step size ($k=0.0125$) (Table~\ref{tab:convergenceEx2}). The accuracy of the proposed scheme is equivalent to the existing fourth-order scheme with reduced computational time being the main source of efficiency (Fig~\ref{fig:efficiencyEx2}C, D). Interestingly, the speed of ETDRK4P22-IF is comparable to the second-order dimensional splitting scheme ETDRDP-IF, but its errors are three orders of magnitude smaller for each time step (Fig~\ref{fig:efficiencyEx2}D, Table~\ref{tab:AppEx2RDPSBDF4}). Similarly, SBDF4 is much slower and less accurate than the proposed scheme at all time steps (Fig~\ref{fig:efficiencyEx2}D, Table~\ref{tab:AppEx2RDPSBDF4}). 
 
\subsection{Enzyme Kinetics}
\label{Ex:Model3}
\begin{figure}[ht]
   \centering
    \includegraphics[width=0.95\textwidth]{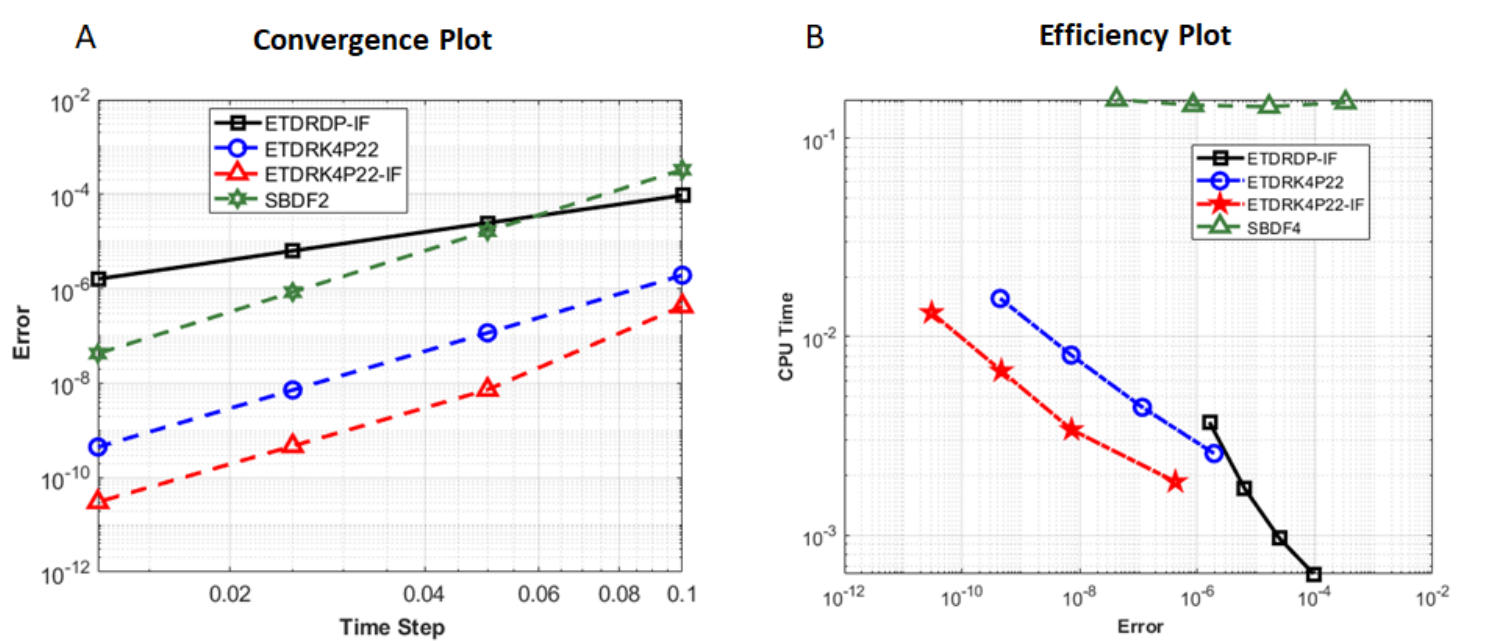}       
 \caption{\textbf{Convergence and efficiency plots for Example~\ref{Ex:Model3}}. Log-log plots showing A) convergence and B) efficiency of ETDRK4P22-IF for the enzyme kinetics problem. Errors are generated between successive numerical solutions with varying temporal step size on a fixed spatial mesh size of $h=0.05$. The diffusion coefficient is $d=0.25$.}
 \label{fig:efficiencyEx3}
 \end{figure}

We consider the two-dimensional analogue of the one-dimensional reaction-diffusion equation with Michaelis-Menten enzyme kinetics \cite{muller2015methods, Bhatt2015,cherruault1990stability}
\begin{equation}
\frac{\partial u}{\partial t} = d \Delta u - \frac{u}{(1+u)}\,,\quad 0<x\,,\,y<1\,,\quad t>0
\end{equation}
with homogeneous Dirichlet boundary conditions and initial condition 
\[ u(x,y,0)=\sin(\pi x)\sin(\pi y)\,,\quad 0\leq x\,,\,y\leq 1\,. \] 
We use this problem to evaluate the performance of the ETDRK4P22-IF scheme for non-linear scalar RDEs with Dirichlet boundary conditions and compare its performance to the existing ETDRK4P22 scheme and ETDRDP-IF schemes. The spatial derivative is discretized as described in Example \ref{Ex:Model3}. We note that this problem has no exact solution and so all errors shown are course-to-fine grid errors as described in the introduction of this section. Our empirical convergence studies show that the scheme is fourth-order accurate with better accuracy than the existing fourth-order ETDRK4P22 scheme (Table~\ref{tab:convergenceEx3}, Fig~\ref{fig:efficiencyEx3}A) and the second-order ETDRDP-IF scheme (Fig.~\ref{fig:efficiencyEx3}A, Table~\ref{tab:AppEx3RDPSBDF4}). Whereas the new scheme is more computational efficient (Table~\ref{tab:convergenceEx2}, Fig.~\ref{fig:efficiencyEx3}B) than the existing ETDRK4P22 scheme, the savings in computational time is not as substantial as was observed in the previous problem. This is most likely because our convergence analysis is performed at a fixed spatial mesh size of $h=0.05$. At this resolution, the computational cost in solving the linear systems in the scheme is relatively small, thus providing little advantage for the dimensional splitting. Whereas ETDRK4P22-IF is generally more accurate than ETDRK4P22, for relatively small problems (less than $20\times 20$ system), it provides little to no speed up in CPU time at this scale. Observe that for all time steps SBDF4 has a substantially high computational cost compared to the other schemes. In fact, there appears to be very little difference in the CPU time as we decrease the time step (Fig.~\ref{fig:efficiencyEx3}B, Table~\ref{tab:AppEx3RDPSBDF4}). This is due to the very small time step needed to compute the starting values of the scheme, as pointed out in the previous examples. Thus, the one-step nature of fourth order ETD schemes make them more computationally efficient than fourth-order IMEX schemes. 

 \begin{table}[ht]
 \begin{center}
      \caption{Convergence table showing fourth-order convergence of the new ETDRK4P22-IF compared with the existing ETDRK4P22 scheme for the enzyme kinetics problem.}
 \resizebox{\textwidth}{!}{
 \begin{tabular}{|c|c|c|c|c|c|c|c|}
 \hline
      &  & \multicolumn{3}{|c|}{ETDRK4P22-IF} & \multicolumn{3}{|c|}{ETDRK4P22} \\
  $k$ & $h$ & Error& $\tilde{p}$ & Cpu Time (sec) & Error & $\tilde{p}$ & Cpu Time (sec) \\

 \hline
 0.10 & 0.05 & 4.2433$\times 10^{-7}$ &--&0.003& 1.9274$\times 10^{-6}$ & -- &0.004\\
 \hline
 0.05  & 0.05 & 7.2737 $\times 10^{-9}$ &5.87& 0.007 & 1.1628$\times 10^{-7}$ & 4.05 &0.008\\
\hline
 0.025  & 0.05 & 4.666$\times 10^{-10}$ &3.96& 0.012 &7.1638$\times 10^{-9}$ & 4.02 &0.014\\
\hline
 0.0125  & 0.05 & 3.0407$\times 10^{-11}$ &3.94& 0.023& 4.4488$\times 10^{-10}$ & 4.01 &0.028\\
 \hline
  \end{tabular}
  }
  \end{center}
    \label{tab:convergenceEx3}
  \end{table}

 \subsection{Non-smooth enzyme kinetics with mismatched initial and boundary data}
\label{Ex:Model4}
Here, we evaluate the performance of ETDRK4P22-IF for solving non-smooth problems by presmoothing with the lower order scheme presented in Section \ref{finalscheme}. Such pre-smoothing techniques have been used to improve the damping properties of the Crank-Nicolson scheme \cite{wade2007smoothing}. The test problem is similar to the model presented in Example \ref{Ex:Model3} with a diffusion coefficient $d=1$ and the following initial condition

\[ u(x,y,0)= 1. \] 
This initial condition is discontinuous at the boundary where homogeneous Dirichlet boundary conditions are applied. Such mismatch of initial and boundary data are known to generate spurious oscillations which, if not properly damped at the initial stages of the simulation, can significantly affect the accuracy of solutions. 

\begin{figure}[h]
   \centering
    \includegraphics[width=0.95\textwidth]{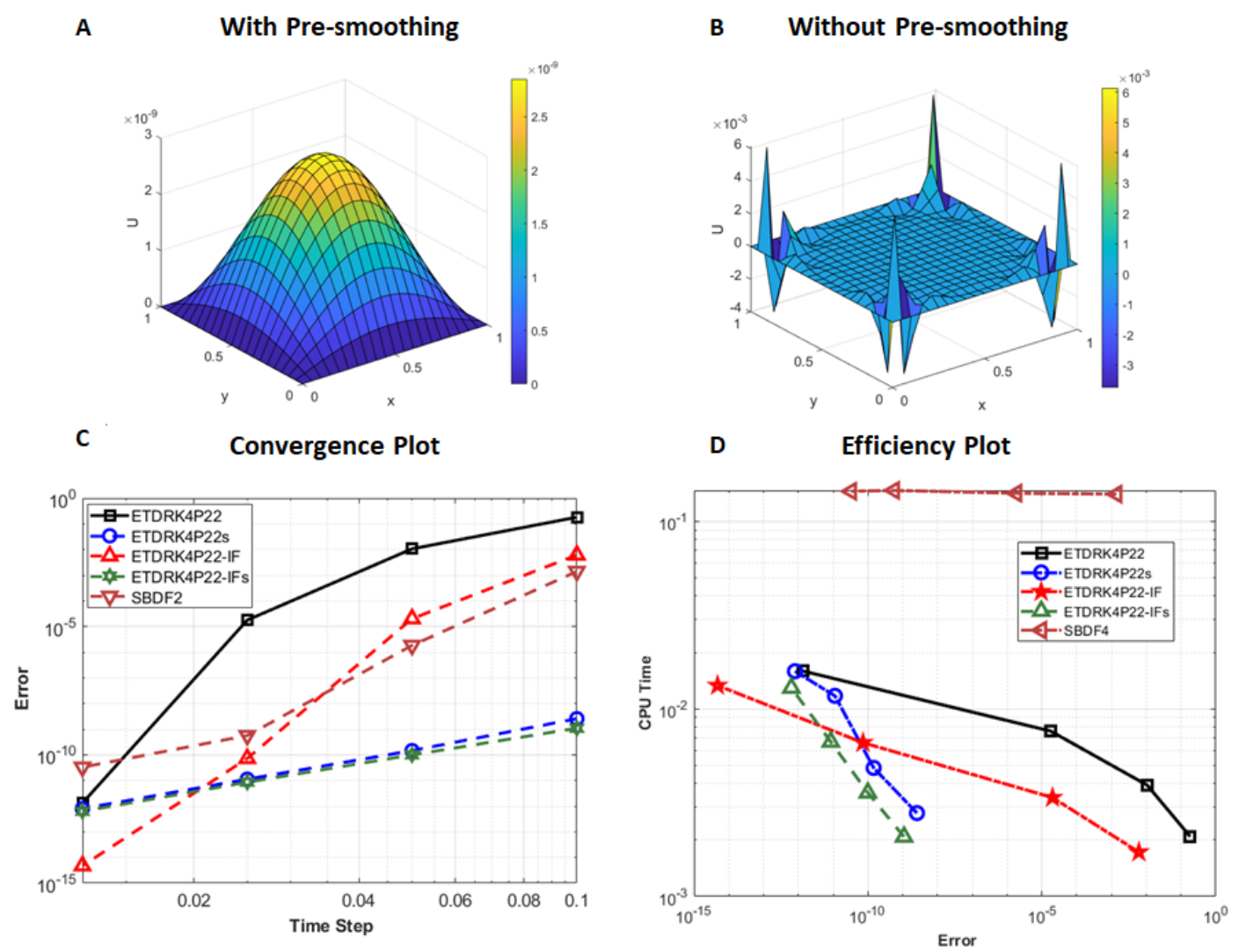}       
 \caption{\textbf{Numerical solution for Example~\ref{Ex:Model3} with mismatched initial and boundary data}. A) Solution profile using ETDRK4P22-IF with three steps of presmoothing. B) Solution profile using ETDRK4P22-IF without presmoothing.  C) Convergence Plots. D) Efficiency Plots. All numerical solutions are generated on a spatial grid with mesh size $h=0.05$. A time step of  $k=0.1$ was used for the solutions profiles in A and B. The diffusion coefficient for the model is taken as $d=1$ in all cases.}
 \label{fig:solutionEx4}
 \end{figure}

 \begin{table}[h]
 \begin{center}
  \caption{Convergence of ETDRK4P22-IF on the nonsmooth enzyme kinetics problem with and without initial smoothing steps. ETDRK4P22-IFs denotes the scheme with three presmoothing steps}
 \resizebox{\textwidth}{!}{
 \begin{tabular}{|c|c|c|c|c|c|c|c|}
 \hline
      &  & \multicolumn{3}{|c|}{ETDRK4P22-IF} & \multicolumn{3}{|c|}{ETDRK4P22-IFs} \\
  $k$ & $h$ & Error& $\tilde{p}$ & Cpu Time (sec) & Error & $\tilde{p}$ & Cpu Time (sec)\\

 \hline
 0.10 & 0.05 & 6.1306$\times 10^{-3}$ &-- &0.002& 1.0894$\times 10^{-9}$ & -- &0.002\\
 \hline
 0.05  & 0.05 & 2.0160 $\times 10^{-5}$ &8.25& 0.003 & 9.9321$\times 10^{-11}$ & 3.46 &0.004\\
\hline
 0.025  & 0.05 & 7.2147$\times 10^{-11}$ &18.09& 0.006 &8.5536$\times 10^{-12}$ & 3.54 &0.007\\
\hline
 0.0125  & 0.05 & 4.7483$\times 10^{-15}$ &13.89& 0.013& 6.2814$\times 10^{-13}$ & 3.77 &0.013\\
 \hline
  \end{tabular}
  }
  \end{center}
    \label{tab:convergenceEx4}
  \end{table}
Without smoothing, the ETDRK4P22-IF scheme fails to damp out spurious oscillations for relatively large time steps (Fig.~\ref{fig:solutionEx4}B, $k=0.1$) and struggles to maintain fourth-order accuracy (Table~\ref{tab:convergenceEx4}). However, when a sufficiently small time step is used ($k=0.0125$), the scheme is able to recover the expected accuracy (Fig.~\ref{fig:solutionEx4}C,D). Performing just three smoothing steps \cite{wade2007smoothing} is sufficient to recover the desired solution profile at the relatively course time step $k=0.1$ (Fig.~\ref{fig:solutionEx4}A) and preserves the fourth-order accuracy of the scheme (Table~\ref{tab:convergenceEx4}). The dimensionally split scheme with smoothing (ETDRK4P22-IFs) is also more efficient than the unsplit scheme (ETDRK4P22) with or without presmoothing for all time steps examined (Fig.~\ref{fig:solutionEx4}D). Again, the SBDF4 scheme is less efficient than all the ETD schemes considered  (Fig.~\ref{fig:solutionEx4}B, Table~\ref{tab:AppEx4RDPSBDF4}). 

% Our empirical convergence studies show that the scheme is fourth-order accurate with better accuracy than the existing fourth-order ETDRK4P22 scheme (Table~\ref{tab:convergenceEx3}, Fig~\ref{fig:efficiencyEx3}A) and second-order ETDRDP-IF scheme (Fig.~\ref{fig:efficiencyEx3}A). The new scheme is also more computational efficient (Table~\ref{tab:convergenceEx2}, Fig.~\ref{fig:efficiencyEx3}B), however the savings in computational time is not as substantial as was observed in the previous problem because our convergence analysis was performed at a fixed spatial mesh size of $h=0.05$. At this resolution, the computational cost in solving the linear systems in the scheme is relatively small, thus providing little advantage for the dimensional splitting. Indeed at a spatial resolution of $h=0.08$ in Example~\ref{Ex:Model1}, ETDRK4P22-IF is about three times faster than ETDRK4P22. Whereas ETDRK4P22-IF is generally more accurate than ETDRK4P22, for relatively small problems (less than $20\times 20$ system), it provides little to no speed up in CPU time. 

 %  \begin{figure}[ht]
 %   \centering
 %  %  \includegraphics[width=0.95\textwidth]{Figures/Figure4.png}       
 % \caption{\textbf{Convergence and efficiency plots for Example~\ref{Ex:Model3}}. Log-log plots showing A) convergence and B) efficiency of ETDRK4P22-IF with and without presmoothing for the non-smooth enzyme Kinetics problem. Errors are generated on a fixed spatial mesh size of $h=0.05$ between successive numerical solutions with varying temporal step size. The diffusion coefficient is $d=0.25$.}
 % \label{fig:efficiencyEx4}
 % \end{figure}

\subsection{Brusselator in 2D}
\label{Ex:Model5}
Finally, we examine the performance of ETDRK4P22-IF for simulating a system of non-linear reaction-diffusion equations, with the Brusselator model \cite{Zegeling2004} as a test case. The model is given by
\begin{align*}
\frac{\partial u}{\partial t} &= \epsilon_1 \Delta u + A + u^2v-(B+1)u\,,\\
\frac{\partial v}{\partial t} &= \epsilon_2 \Delta v + Bu - u^2v\,,
\end{align*}
 in 2D, where the diffusion coefficients are $ \epsilon_1=\epsilon_2=2\cdotp 10^{-3}$, and the chemical parameters $A=1$, $B=3.4$. At the boundary of the domain homogeneous Neumann conditions are imposed. 
%Is this really needed here? Just write homogeneous Neumann.
% \[ \frac{\partial u}{\partial n} |_{\partial \Omega} = \frac{\partial v}{\partial n} |_{\partial \Omega}=0 \hspace{5mm} \text{ with } \Omega = [0,1]^2,\ t\in [0,T]\]

\begin{figure}[h]
   \centering
    \includegraphics[width=0.95\textwidth]{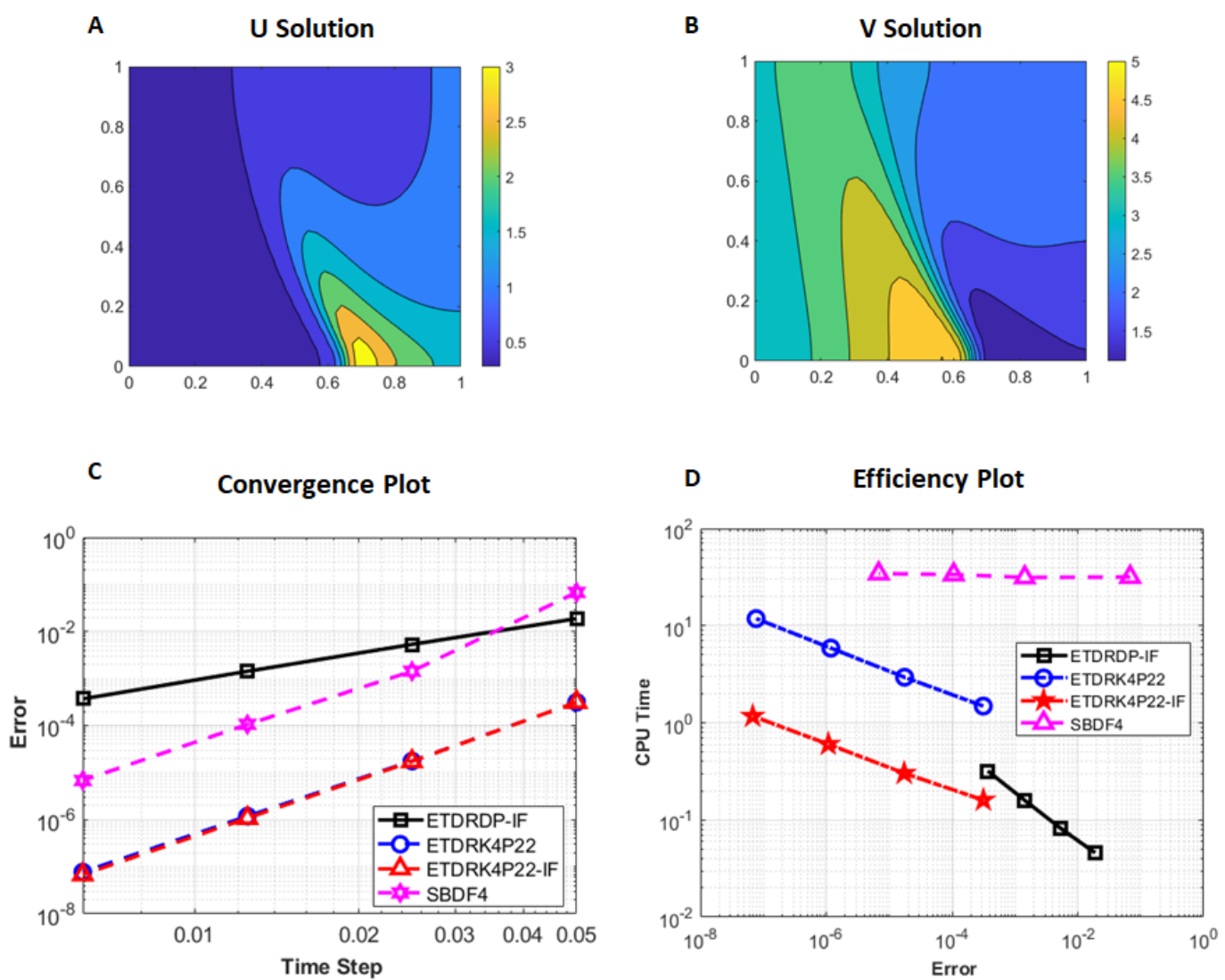}       
 \caption{\textbf{Convergence and efficiency plots for Example~\ref{Ex:Model5}}. A) Numerical solution $U$, B) Numerical solution $V$, log-log plots showing C) Fourth-order convergence, and D) Efficiency of the new fourth-order ETDRK4P22-IF scheme compared with the existing fourth-order ETDRK4P22, second-order ETDRDP-IF and fourth order SBDF4 schemes. Numerical solutions are generated using $40$ nodes per spatial dimension and time step of $0.05$ with final time $T=2$. }
 \label{fig:efficiencyEx5}
 \end{figure}

The initial conditions are given by
\[ u(x,y,0)=\frac{1}{2} + y\,,\quad v(x,y,0)=1+5x\,.\]
Spatial derivatives are discretized as described in Example \ref{Ex:Model2}. Again, this problem has no exact solution and so all errors shown are course-to-fine grid errors as described in the introduction of this section. The numerical solution for each component is shown in Fig~\ref{fig:efficiencyEx5}A,B. Solution profiles are similar to those produced in \cite{Zegeling2004}. Our empirical convergence studies show that the scheme is fourth-order accurate with similar accuracy to the existing fourth-order ETDRK4P22 scheme (Table~\ref{tab:convergenceEx5}, Fig~\ref{fig:efficiencyEx5}C). However, the new scheme is more computationally efficient than ETDRK4P22, with seven times speed up in CPU time at the spatial mesh size of $h=0.0125$ and temporal step size of $k=0.00625$ (Table~\ref{tab:convergenceEx5}, Fig.~\ref{fig:efficiencyEx5}D). As expected, ETDRK4P22-IF is more accurate than the existing second-order ETDRDP-IF scheme, providing errors that are about four orders of magnitude smaller (Fig.~\ref{fig:efficiencyEx5}B, Table~\ref{tab:AppEx5RDPSBDF4}). ETDRK4P22-IF is also more efficient than SBDF4 (Fig.~\ref{fig:efficiencyEx5}B, Table~\ref{tab:AppEx5RDPSBDF4}) .

 \begin{table}[h]
 \begin{center}
   \caption{Convergence table showing fourth-order convergence of the new ETDRK4P22-IF scheme compared with the existing ETDRK4P22 scheme for Brusselator Example~\ref{Ex:Model5} with $T=2$.}
 \resizebox{\textwidth}{!}{
 \begin{tabular}{|c|c|c|c|c|c|c|c|}
 \hline
      &  & \multicolumn{3}{|c|}{ETDRK4P22-IF} & \multicolumn{3}{|c|}{ETDRK4P22} \\
  $k$ & $h$ & Error& $\tilde{p}$ & Cpu Time (sec) & Error & $\tilde{p}$ & Cpu Time (sec) \\

 \hline
 0.05& 0.0125 & 3.1532$\times 10^{-4}$ &-- &0.160& 3.1384$\times 10^{-4}$ & -- &1.486\\
 \hline
 $0.025$  & 0.0125 & 1.7359$\times 10^{-5}$ &4.18& 0.303 & 1.7592$\times 10^{-5}$ & 4.16 & 2.948\\
\hline
 $0.0125$  & 0.0125 & 1.0814 $\times 10^{-6}$ &4.00& 0.603 &1.1859 $\times 10^{-6}$ & 3.89 & 5.879\\
\hline
 $0.00625$  & 0.0125 & 6.7987$\times 10^{-8}$ &3.99 & 1.172 & 7.6968 $\times 10^{-8}$ & 3.95 & 11.836\\
 \hline
  \end{tabular}
  }
  \end{center}
    \label{tab:convergenceEx5}
  \end{table}

\section{Conclusions} 
\label{sec:conclusion}
In this article, we have developed a fourth-order ETD scheme with dimensional splitting that uses a Pad\'e (2,2) rational function to approximate the matrix exponential. The order of convergence of the scheme has been validated using two linear problems with exact solution and two non-linear problems with Dirichlet and Neumann boundary conditions in 2D. In all cases the new scheme is more computational efficient than competing schemes. We expect the same savings to persist for problems with periodic boundary conditions with additional savings introduced through parallelization of the new scheme. Addition of up to three pre-smoothing steps using a third-order L-stable scheme effectively damps out spurious oscillations arising from mismatched initial and boundary data and preserves the fourth-order accuracy of the scheme. Future work will examine 3D applications and investigate the assumptions on the differential operator and non-linear source term that ensure fourth-order accuracy and positivity of the numerical solution.
 \section*{Declaration of Interest Statement}
None.
  \section*{Author Contributions Statement}
\textbf{E.O. Asante-Asamani}: Conception and design, Analysis and interpretation of data, Writing original draft; \textbf{B.A. Wade:} Analysis and interpretation of data, Writing-review \& editing. \textbf{A. Kleefeld:} Analysis and interpretation of data, Writing-review \& editing.
\clearpage

\clearpage
% \bibliographystyle{elsarticle-num}
%\bibliography{ETDRK4}

%%%%%%%%%%%%%%%%%%%%%%%%%%%%%%%%%%%%%%%%%%%%%%%%%

\section*{Appendices}
\label{sec:appendix}

\begin{appendices}
      \setcounter{table}{0}
\renewcommand{\thetable}{A.\arabic{table}}
    \section{Supplementary tables}

 \begin{table}[h]
 \begin{center}
   \caption{Convergence table showing fourth-order accuracy of SBDF4 scheme and second-order accuracy of ETDRDP-IF for Example~\ref{Ex:Model1}}
   \label{tab:AppEx1RDPSBDF4}
 \resizebox{\textwidth}{!}{
 \begin{tabular}{|c|c|c|c|c|c|c|c|}
 \hline
      &  & \multicolumn{3}{|c|}{ETDRDP-IF} & \multicolumn{3}{|c|}{SBDF4} \\
  $k$ & $h$ & Error& $p$ & Cpu Time (sec) & Error & $p$ & Cpu Time (sec) \\

 \hline
  0.1000& 0.0785 & 3.3043$\times 10^{-4}$ &-- &0.008& 2.2150$\times 10^{-4}$ &-- &1.145\\
 \hline
  0.0500 & 0.0393 & 8.4134$\times 10^{-5}$ &1.97& 0.018 & 1.2419$\times 10^{-5}$ &4.16& 14.823\\
\hline
 0.0250  & 0.0196 & 2.1234 $\times 10^{-5}$ &1.99& 0.127& 7.752 $\times 10^{-7}$ &4.00& 115.230\\
\hline
0.0125  & 0.0098 & 5.3332$\times 10^{-6}$ &1.99 & 1.167 & 6.1782$\times 10^{-8}$ &3.65 & 1209.314\\
 \hline
  \end{tabular}
  }
  \end{center}
  \end{table}

\begin{table}[h]
 \begin{center}
   \caption{Convergence table showing fourth-order accuracy of SBDF4 scheme and second-order accuracy of ETDRDP-IF for Example~\ref{Ex:Model2}}
       \label{tab:AppEx2RDPSBDF4}

 \resizebox{\textwidth}{!}{
 \begin{tabular}{|c|c|c|c|c|c|c|c|}
 \hline
      &  & \multicolumn{3}{|c|}{ETDRDP-IF} & \multicolumn{3}{|c|}{SBDF4} \\
  $k$ & $h$ & Error& $p$ & Cpu Time (sec) & Error & $p$ & Cpu Time (sec) \\

 \hline
  0.1000& 0.1496 & 4.4405$\times 10^{-4}$ &-- &0.002& 2.3248$\times 10^{-4}$ &-- &0.189\\
 \hline
  0.0500 & 0.0766 & 1.0839$\times 10^{-4}$ &2.03& 0.010 & 1.3094$\times 10^{-5}$ &4.15& 1.198\\
\hline
 0.0250  & 0.0388 & 2.6811 $\times 10^{-5}$ &2.02& 0.073& 8.1722 $\times 10^{-7}$ &4.00& 14.316\\
\hline
0.0125  & 0.0195 & 6.6694$\times 10^{-6}$ &2.01 & 1.809 & 6.4410$\times 10^{-8}$ &3.67 & 114.896\\
 \hline
  \end{tabular}
  }
  \end{center}
  \end{table}
  
\begin{table}[h]
 \begin{center}
   \caption{Convergence table showing fourth-order accuracy of SBDF4 scheme and second-order accuracy of ETDRDP-IF for Example~\ref{Ex:Model3}}
       \label{tab:AppEx3RDPSBDF4}
 \resizebox{\textwidth}{!}{
 \begin{tabular}{|c|c|c|c|c|c|c|c|}
 \hline
      &  & \multicolumn{3}{|c|}{ETDRDP-IF} & \multicolumn{3}{|c|}{SBDF4} \\
  $k$ & $h$ & Error& $\tilde{p}$ & Cpu Time (sec) & Error & $\tilde{p}$ & Cpu Time (sec) \\

 \hline
  0.1000& 0.1653 & 9.4569$\times 10^{-5}$ &-- &0.0006& 3.3077$\times 10^{-4}$ &-- &0.1511\\
 \hline
  0.0500 & 0.1653 & 2.4719$\times 10^{-5}$ &1.94& 0.0009 & 1.6554$\times 10^{-5}$ &4.32& 0.1437\\
\hline
 0.0250  & 0.1653 & 6.3326$\times 10^{-6}$ &1.96& 0.0017& 8.4646 $\times 10^{-7}$ &4.29& 0.1464\\
\hline
0.0125  & 0.1653 & 1.6036$\times 10^{-6}$ &1.98 & 0.0037 & 4.2354$\times 10^{-8}$ &4.32 & 0.1560\\
 \hline
  \end{tabular}
  }
  \end{center}
  \end{table}

  \begin{table}[h]
 \begin{center}
   \caption{Convergence table for ETDRK4P22 and ETDRK4P22s schemes applied to Example~\ref{Ex:Model4}}
       \label{tab:AppEx4RDPSBDF4}
 \resizebox{\textwidth}{!}{
 \begin{tabular}{|c|c|c|c|c|c|c|c|}
 \hline
      &  & \multicolumn{3}{|c|}{ETDRK4P22} & \multicolumn{3}{|c|}{ETDRK4P22s} \\
  $k$ & $h$ & Error& $\tilde{p}$ & Cpu Time (sec) & Error & $\tilde{p}$ & Cpu Time (sec)\\

 \hline
  0.1000& 0.1653 & 1.8184$\times 10^{-1}$ &-- &0.002& 2.5622$\times 10^{-9}$ &-- &0.003\\
 \hline
  0.0500 & 0.1653 & 1.0646$\times 10^{-2}$ &2.09& 0.004 & 1.4589$\times 10^{-10}$ &4.13& 0.005\\
\hline
 0.0250  & 0.1653 & 1.7765$\times 10^{-5}$ &9.23& 0.008& 1.124 $\times 10^{-11}$ &3.70& 0.012\\
\hline
0.0125  & 0.1653 & 1.3239$\times 10^{-12}$ &23.68 & 0.016& 7.9819$\times 10^{-13}$ & 3.82 & 0.016\\
 \hline
  \end{tabular}
  }
  \end{center}
  \end{table}

  \begin{table}[h]
 \begin{center}
   \caption{Convergence table showing fourth-order accuracy of SBDF4 scheme and second-order accuracy of ETDRDP-IF for Example~\ref{Ex:Model5}}
       \label{tab:AppEx5RDPSBDF4}
 \resizebox{\textwidth}{!}{
 \begin{tabular}{|c|c|c|c|c|c|c|c|}
 \hline
      &  & \multicolumn{3}{|c|}{ETDRDP-IF} & \multicolumn{3}{|c|}{SBDF4} \\
  $k$ & $h$ & Error& $\tilde{p}$ & Cpu Time (sec) & Error & $\tilde{p}$ & Cpu Time (sec) \\

 \hline
  0.0500& 0.0125 & 1.8794$\times 10^{-2}$ &-- &0.046& 6.7785$\times 10^{-2}$ &-- &31.672\\
 \hline
  0.0250 & 0.0125 & 5.3233$\times 10^{-3}$ &1.82& 0.082 & 1.4339$\times 10^{-3}$ &5.56& 31.391\\
\hline
 0.0125 & 0.0125 & 1.4270$\times 10^{-3}$ &1.90& 0.160& 1.0627 $\times 10^{-4}$ &3.75& 33.662\\
\hline
0.0063  & 0.0125 & 3.7001$\times 10^{-4}$ &1.95 & 0.310& 6.8328$\times 10^{-6}$ & 3.96 & 34.430\\
 \hline
  \end{tabular}
  }
  \end{center}
  \end{table}
%%%%%%%%%%%%%%%%%%%%%%%%%%%%%%%%%%%%%%%%%%%%%%%%%%%%%%%%%%%%%%%%%%%%%%

\end{appendices}

\end{document}